\documentclass[reqno,11pt]{amsart}
\usepackage{amsthm}
\usepackage{amsmath,mathrsfs}
\usepackage{color}
\usepackage{amssymb}
\usepackage{hyperref}
\usepackage{enumerate}
\usepackage{latexsym,array}
\usepackage{amsfonts}
\usepackage{shadow}
\usepackage[utf8]{inputenc}

\newtheorem{Pa}{Paper}[section]
\newtheorem{Tm}[Pa]{{\bf Theorem}}
\newtheorem{La}[Pa]{{\bf Lemma}}
\newtheorem{Dn}[Pa]{{\bf Definition}}
\newtheorem{Cy}[Pa]{{\bf Corollary}}

\newtheorem{Rk}[Pa]{{\bf Remark}}
\newtheorem{Pn}[Pa]{{\bf Proposition}}

\def\C{\mathbb C}
\def\hh{\mathbb{H}}

\def\loc{\mathrm{loc}}

\author[F. Colombo]{Fabrizio Colombo}
\address{(FC) Politecnico di
Milano\\Dipartimento di Matematica\\Via E. Bonardi, 9\\20133
Milano, Italy}
\email{fabrizio.colombo@polimi.it}

\author[S. Mongodi]{Samuele Mongodi}
\address{(SM) Politecnico di
Milano\\Dipartimento di Matematica\\Via E. Bonardi, 9\\20133
Milano, Italy}
 \email{samuele.mongodi@polimi.it}

\title[]
{The Cauchy Transform in the  slice hyperholomorphic setting and related topics} \oddsidemargin 0.2in \evensidemargin
0.2in \topmargin -0.5in \textwidth 15.5truecm \textheight 23truecm
\def\H{\mathbb H}

\def\R{\mathbb R}
\def\N{\mathbb N}
\def\C{\mathbb C}

\def\(s){\mathscr S(\R\times\R)}



\begin{document}
\maketitle \tableofcontents
\parindent 0cm
\begin{abstract}
In this paper we study the additive splitting associated to the quaternionic Cauchy transform defined by the
Cauchy formula of slice hyperholomorphic functions. Moreover, we introduce and study the analogue of the fundamental solution of the
global operator of slice hyperholomorphic functions. We state our results in the quaternionic setting but several results hold
for Clifford algebra-valued function with minor changes in the proofs.
\end{abstract}

\parindent 0cm

\section{Introduction}
\setcounter{equation}{1}

The classical function theory for quaternionic-valued functions is based on
the so called Cauchy-Fueter operator, see \cite{fueter 1}.
In the case of Clifford algebra-valued functions, the theory of nullsolutions of the Dirac operator is a widely studied theory, see e.g. the books \cite{dss,csss,ghs}.  These functions are called monogenic and since the Dirac operator factorizes the Laplace operator in $n$ dimension they are harmonic.
\\
The more recent theory of slice hyperholomorphic functions has been introduced over ten years ago, see \cite{GS2} and  \cite{slicecss}, and it is an alternative function
theory. These functions are often called slice regular when they are quaternion-valued
while when they are Clifford algebra-valued they are called slice monogenic.
The developments of this theory  can be found in
 in the books \cite{ACS_book,css_book,CGK_book,ENTIRE_book,GSS_book}.
One of the most important applications of slice hyperholomorphic functions is in operator theory, where they play the role of holomorphic functions for the classical spectral theory.
Most of the classical results in operator theory can be extended to quaternionic operators and to
$n$-tuples of non commuting operators using the notion of $S$-spectrum, see \cite{ACS_book,CGK_book,css_book}.
 \\
It is interesting to point out that even though the theory of monogenic and slice hyperholomorphic functions are quite different there are some interesting relations between them which are given by the
Fueter-Sce-Qian mapping theorem, see \cite{CoSaSo1}, or by the Radon and the dual Radon transform,  see  \cite{Radon}.
\\
Many results that hold for scalar valued slice hyperholomorphic functions can be extended to the vector-valued case, introduced in the paper
\cite{ACLS} and extensively studied in \cite{ACS_book}.

In this paper we study continuous functions on a closed boundary of a suitable  open set of the space of
quaternions that admit an additive splitting as sum of two functions,
one continuous and  slice hyperholomorphic inside or outside the open set, respectively.
In the classical case, such results are contained for example in the book \cite{GL}.
As in the complex case, not all continuous functions admit such a splitting. The splitting can be assured when we consider H\"older continuous functions.
In the slice hyperholomorphic case, the multiplicative splitting seems to be unnatural because the composition of slice hyperholomorphic functions is not
always a function of the same type.

 To introduce the splitting, we consider $\partial U $ to be the smooth boundary of a bounded set $U\subset \mathbb{H}$, where $\mathbb H$ is the set of quaternions. Let $f:\partial U \to \mathcal{X}$ be a slice continuous function and $\mathcal{X}$ be a quaternionic Banach space. Let $U_+$, $U_-$ be the inside and the outside of $U$, respectively.
We say that $f$ splits additively with respect to $\partial U$ if there exist functions
$f_{-}:\overline{U}_{-}\cup \{\infty\}\to \mathcal{X}$ and $f_{+}:\overline{U}_{+}\to \mathcal{X}$ where
$f_{-}$ is continuous on $\overline{U}_{-}$ and slice hyperholomorphic in ${U}_{-}\bigcup\{\infty\}$ and $f_+$ is continuous on $\overline{U}_+$ and slice hyperholomorphic in
$U_+$ such that, on $\partial U$ we have
$$
f=f_-+f_+.
$$
In the sequel, $\mathbb{S}$ is the unite sphere of purely imaginary quaternions,  ${\rm Re}(s)$ denotes
the real part of a quaternion $s$ and  $|s|^2$ is the Euclidean squared norm of $s\in \mathbb{H}$. By $\mathbb{C}_j$ we denote the complex plane with imaginary unit $j\in \mathbb{S}$.
Let $\partial (U\cap \mathbb{C}_j)$, $j\in \mathbb{S}$ be the piecewise $\mathcal{C}^1$-contour  of $U\cap \mathbb{C}_j$ where
$U\subset \mathbb{H}$. Let  $f:\partial U \to \mathcal{X}$ be a left slice continuous function.
We define the left Cauchy integral transform with respect to  $\partial U$ of $f$ as
\begin{equation}
 \hat{f}(p)=\frac{1}{2 \pi}\int_{\partial (U\cap \mathbb{C}_j)} S_L^{-1}(s,p)ds_j f(s), \ \ \ p\in  \mathbb{H}\setminus \partial U,
\end{equation}
where the Cauchy kernel
$$
 S_L^{-1}(s,p):=-(p^2-2p{\rm Re}(s)+|s|^2)^{-1}(p-\overline{s})
$$
is a left  slice hyperholomorphic function in the variable $p$ and right slice hyperholomorphic function in $s$.
\\
Let us mention that there are different ways to introduce slice hyperholomorphic functions, one of which involves a non constant
coefficients differential operator that was introduced in
\cite{GLOBAL}
 \begin{equation}\label{glop}
 G_Lf(q)=|\underline{q}|^2\dfrac{\partial f}{\partial x_0}(q)+\underline{q}\sum_{j=1}^3 x_j\dfrac{\partial f}{\partial x_j}(q).
 \end{equation}
 It is important to study the counterpart of the fundamental solution of this operator in terms of the Cauchy kernel. In this work we show that
$$
G_L(S^{-1}_L(s,\cdot))=2\pi j|\underline{s}|^2\delta_s
$$
in the sense of distribution, where $\delta_s$ is the Dirac delta function.
\\
Finally, we note that in this paper we consider the quaternionic setting but many results may be extended to the Clifford algebra setting with minor changes.
\subsection*{Plan of the paper.}
In Section \ref{due} we recall some basic facts on slice hyperholomorphic functions quaternionic-valued
 that can be naturally extended to vector-valued functions.
In Section \ref{tre} we study the additive splitting for slice continuous functions through Cauchy transform.
The case of H\"older continuous functions is treated in Section \ref{quattro}.
In Section \ref{cinque} we study the fundamental solution of the global operator \eqref{glop} and some related topics.

\section{Preliminary results on slice operator-valued functions}\label{due}
Let $\mathbb{H}$ be the algebra of quaternions.
The imaginary units  in $\mathbb{H}$, denoted by $e_1$, $e_2$, $e_3$, satisfy the relations
$ e_1^2=e_2^2=e_3^2=-1$, $e_1e_2 =-e_2e_1 =e_3$, $e_2e_3 =-e_3e_2 =e_1$, $e_3e_1 =-e_1e_3 =e_2$
 and an element in $\hh$ is of the form $q=x_0+e_1x_1+e_2x_2+e_3x_3$, for $x_\ell\in\mathbb{R}$.
The real part, the imaginary part (or vector part) and the modulus of a quaternion are defined as
$ {\rm Re}(q)=x_0$, ${\rm Im}(q)=e_1x_1+e_2x_2+e_3x_3$,
$|q|^2=x_0^2+x_1^2+x_2^2+x_3^2$, respectively.
The conjugate of the quaternion $q=x_0+e_1x_1+e_2x_2+e_3x_3$ is defined by
$$
\bar q={\rm Re }(q)-{\rm Im }(q)=x_0-e_1x_1-e_2x_2-e_3x_3.
$$
Let us denote by $\mathbb{S}$ the unit sphere of purely imaginary $\mathbb{S}$
quaternions, i.e.
$$
\mathbb{S}=\{q=e_1x_1+e_2x_2+e_3x_3\ {\rm such \ that}\
x_1^2+x_2^2+x_3^2=1\}.
$$
Given a non-real quaternion $q$ we can write it as $q=x_0+{\rm Im} (q)=x_0+j |{\rm Im} (q)|$, where $j={\rm Im} (q)/|{\rm Im} (q)|\in\mathbb{S}$. Furthermore, we can associate to $q$ the 2-dimensional sphere $[q]$ defined by
$$
[q]=\{x_0+j  |{\rm Im} (q)| \ : \ j \in\mathbb{S}\}.
$$

In order to adapt some results known in the case of slice hyperholomorphic functions quaternionic-valued to
 the case of vector-valued slice hyperholomorphic functions we need some definitions and preliminary results.
 \\
  An open set $U \subset\hh$ is said to be axially symmetric if $[q]\subset U$ whenever $q\in U$.
\begin{Dn}
Let $\mathcal{X}$ is a two-sided quaternionic Banach algebra and let $U\subseteq\hh$ be an axially symmetric open set.
 Let $\mathcal{U}\subseteq\mathbb{R}\times \mathbb{R}$ be such that $p=u+j v\in U$ for all $(u,v)\in\mathcal{U}$ and for all $j\in \mathbb{S}$.
Functions of the form
\begin{equation}\label{SliceL}
f(p)=f(u+jv)=f_0(u,v) +jf_1(u,v)
\end{equation}
 where $f_0,f_1:\ \mathcal U \to \mathcal{X}$ depend only on $u,v$, and satisfy
 $$
 f_0(u,-v)=f_0(u,v),\ \ \  f_1(u,-v)=-f_1(u,v)
 $$
 are called slice functions. Moreover:
 \begin{enumerate}
 \item[(i)]
 We say that $f$  is a left slice $\mathcal{C}^k$-function  if  $f_0$ and $f_1$ are of class  $\mathcal{C}^k$, for $k\in \mathbb{N}_0$.
 We will denote this class of functions by $\mathcal{SC}^k_L(U,\mathcal{X})$.
 \item[(ii)]
 If  $f_0$ and $f_1$ are real  differentiable and  satisfy the Cauchy-Riemann
equations
$$
\partial_v f_0(u,v) +\partial_uf_1(u,v)=0,\ \ \   \partial_u f_0(u,v) -\partial_vf_1(u,v)=0,
$$
we will call $f$ left slice hyperholomorphic.
 We will denote this class of functions by
$\mathcal{SH}_L(U,\mathcal{X})$.
\end{enumerate}
 In the case $f$ is of the form
$$
f(p)=f(u+jv)=f_0(u,v) +f_1(u,v)j
$$
then the set of right slice $\mathcal{C}^k$-functions will be denoted by $\mathcal{SC}^k_R(U,\mathcal{X})$ and the set of
right slice hyperholomorphic functions will be denoted by $\mathcal{H}_R(U,\mathcal{X})$.
\end{Dn}

\begin{Tm}[Representation Formula]\label{RepFo}
Let $\mathcal{X}$ be a two-sided quaternionic Banach algebra.
Let $U\subset \mathbb{H}$ be an axially symmetric open set and let $j\in\mathbb{S}$. For any
$p = p_0 + j_p p_1\in U$ set $p_j := p_0 + jp_1$. If $f\in\mathcal{SC}^k_L(U,\mathcal{X})$, $k\in\mathbb N_0$, then
\[
f(p) = \frac12(1-j_pj)f(p_j) + \frac12(1+j_pj)f(\overline{p_j}) \quad\text{for all \ \ $p\in U.$}
\]
If $f\in \mathcal{SC}^k_R(U,\mathcal{X})$, then
\[f(p) = f(p_j)(1-jj_p)\frac12 + f(\overline{p_j})(1+jj_p)\frac12 \quad\text{for all\ \  $p\in U$}.
\]
\end{Tm}
We note that the representation formula  holds for slice functions $f(p)=f(u+jv)=f_0(u,v) +jf_1 (u,v)$ and does not depend
on the regularity of the functions $f_0(u,v)$ and $f_1(u,v)$. It is an easy consequence of the definition, see \cite{GP}.

\begin{Dn}
Let $U \subset\hh$ be an axially symmetric open set and let $\mathcal{X}$ be a two-sided quaternionic Banach algebra.
Let $f,g\in\mathcal{SH}_L(U,\mathcal{X})$ with $f(q) = f_0(u,v) + j f_1(u,v)$ and $g(q) = g_0(u,v) + j g_1(u,v)$ for $q = u + jv\in U$; we define their  left slice hyperholomorphic product as
\begin{equation}\label{starleft}
\begin{split}
f \star_l g := &f_0g_0 - f_1g_1 + j\left(f_0g_1 + f_1g_0\right).
\end{split}
\end{equation}
 Let $f,g\in\mathcal{SH}_R(U,X)$ with $f(q) = f_0(u,v) +  f_1(u,v)j$ and $g(q) = g_0(u,v)  + g_1(u,v)j$ for $q = u + jv\in U$; we define their right slice hyperholomorphic product as
\begin{equation}\label{starright}
\begin{split}
f \star_r g := &f_0g_0 - f_1g_1 + \left(f_0g_1 + f_1g_0\right)j.
\end{split}
\end{equation}
\end{Dn}

We recall that we have the notion of weakly slice hyperholomorphic function also in also in this specific setting. \begin{Dn}
Let $U \subset \hh$ be an axially symmetric open set.
A function $f:U\to \mathcal{X}$ with values in a quaternionic Banach space $\mathcal{X}$ is called
 weakly left slice hyperholomorphic if $\Lambda f$ is left slice hyperholomorphic for any $\Lambda\in \mathcal{X}'$, where $\mathcal{X}'$ is the dual space of $\mathcal{X}$. A similar definition holds for
  weakly right slice hyperholomorphic.
\end{Dn}
The following version of
Liouville's theorem is an easy consequence of the analogue result in the quaternionic valued setting (see \cite{GSS_book}):
\begin{Tm}[Liouville]
Let $\mathcal{X}$ be a Banach space and let $\mathcal{X}'$ be the dual space of $\mathcal{X}$. Suppose that
$f:\hh\to \mathcal{X}$ be a slice left hyperholomorphic function. Suppose that for every $\Lambda\in \mathcal{X}'$ the function
$\Lambda f$ is bounded on  $\hh$. Then $f$ is constant.
\end{Tm}

\begin{Tm}
Let $\mathcal{X}$ be a quaternionic two-sided Banach space,  let $U$ be an open axially symmetric subset
of $\hh$ and let $f: U\to \mathcal{X}$ be a real differentiable left slice function.
Then the function $f$ is strongly left slice hyperholomorphic if and only if
the function $f$ admits left slice derivative, that is
\begin{equation}\label{VecSDeriv}
\partial_S f(q) := \lim_{p\to q, p\in \mathbb{C}_j} (p-q)^{-1}(f(p)-f(q))
\end{equation}
exists for all $q = u + j v\in U$ in the topology of $\mathcal{X}$ and it exists for any $j\in\mathbb{S}$ if $q$ is real.
\end{Tm}
A similar statement holds for  strongly right slice hyperholomorphic functions.
In this case the right slice derivative, is defined by
\[
\partial_S f(q) := \lim_{p\to q, p\in \mathbb{C}_j} (f(p)-f(q))(p-q)^{-1}
\]
exists for all $q = u + j v\in U$ in the topology of $\mathcal{X}$ and it exists for any $j\in\mathbb{S}$ if $q$ is real.

\begin{Tm}[Maximum Principle]
Let $\mathcal X$ be a  two--sided quaternionic Banach space and let
  $\overline{U}\subset \hh$ be an axially symmetric domain.
  Let $f:\overline{U}\to \mathcal{X}$ be a continuous function that is slice left (or right) hyperholomorphic in $U$ with values in $\mathcal{X}$.
  Then
  $$
  \max_{p\in \overline{U}}\|f(p)\|=\max_{p\in \partial U}\|f(p)\|.
  $$
\end{Tm}
\begin{proof}
As in the classical case it follows from the  Maximum Principle of slice hyperholomorphic functions with values in $\hh$, see \cite{GSS_book} and the quaternionic version of the Hahn-Banach theorem, see \cite{css_book}.
\end{proof}

As a consequence one deduces:
\begin{Tm}\label{1.3.6}
Let $\partial(U\cap \mathbb{C}_j)$, $j\in \mathbb{S}$ be a $\mathcal{C}^1$-contour in $\mathbb{C}_j$ and let $f$ be a right slice continuous
function $f:\partial(U\cap \mathbb{C}_j)\to \mathcal{X}$ where $\mathcal{X}$ is a quaternionic Banach space. Then
$$
\left\|\int_{\partial(U\cap \mathbb{C}_j)} f(s)ds_j \right\|\leq |\partial(U\cap \mathbb{C}_j)| \max_{s\in \partial(U\cap \mathbb{C}_j)}\|f(s)\|,
$$
where $|\partial(U\cap \mathbb{C}_j)|$ denotes the lenght of $\partial(U\cap \mathbb{C}_j)$.
\end{Tm}
 The proof of this result closely follows the one in the complex case, see Proposition 1.3.6 in \cite{GL} and thus will be omitted.

Vector-valued slice hyperholomorphic functions admit power series expansions at real point that are similar to those for classical complex holomorphic functions. Specifically:
\begin{Tm}
Let $\mathcal X$ be a two sided quaternionic Banach space.
Let $\alpha\in\mathbb{R}$, let $r>0$ and let $B_{r}(\alpha) = \{q\in\hh: |q-\alpha|<r\}$. If $f\in\mathcal{SH}_L(B_r(\alpha), \mathcal{X})$ then
\begin{equation}
\label{PowSerL}
f(q) = \sum_{n= 0}^{+\infty} (q-\alpha)^n\frac{1}{n!}\partial_S^n f(\alpha)\qquad \forall q = u + j v \in B_r(\alpha).
\end{equation}
If $f\in\mathcal{SH}_L(B_r(\alpha), \mathcal{X})$, then
\[
f(q) = \sum_{n= 0}^{+\infty}\frac{1}{n!}\left(\partial_S^n f(\alpha)\right)  (q-\alpha)^n\qquad \forall q = u + j v \in B_r(\alpha).
\]
\end{Tm}
We now discuss how to generalize the theory of Laurent series to slice hyperholomorphic functions with values in a two sided quaternionic Banach space $\mathcal{X}$.
Given the series
 $\sum_{n \in \mathbb{Z}} q^n f_n$, for $f_n\in \mathcal{X}$ we denote by $A(0, R_1, R_2)$ the four-dimensional spherical shell
$$
A(0, R_1, R_2) = \{q \in \mathbb{H} : R_1 < |q| < R_2\}.
$$
We have
\begin{Tm}\label{laurentregularity}
Let $\sum_{n \in \mathbb{Z}} q^n f_n$, $f_n\in\mathcal X$ be a series having domain of convergence
$A = A(0,R_1,R_2)$ with $R_1 < R_2$. Then
$f : A \to \mathbb{H} \ \ q \mapsto \sum_{n \in \mathbb{Z}} q^n f_n$ is a left slice hyperholomorphic function.
\end{Tm}
\begin{proof}
It follows from the fact that
if $\{f_n\}_{n \in \mathbb{Z}} \subset \mathcal{X}$, there exist $R_1, R_2$ with $0
\leq R_1< R_2 \leq \infty$ such that
  the series  $\sum_{n \in \mathbb{N}} q^n f_n$ and
$ \sum_{n \in \mathbb{N}} q^{-n} a_{-n}$ both converge
absolutely and uniformly on compact subsets of $A = A(0,R_1,R_2)$;
while for all $q \in \mathbb{H} \setminus \bar A$,
either $\sum_{n \in \mathbb{N}} q^n f_n$ or $ \sum_{n \in \mathbb{N}} q^{-n} a_{-n}$ is divergent. Where the series converges it is trivially slice hyperholomorhic.
\end{proof}
Conversely, we have that all left slice hyperholomorphic  functions $f : A(0,R_1,R_2) \to \mathcal{X}$ admit
Laurent series expansions.
\begin{Tm}[Laurent Series Expansion]\label{expansion}
Let  $A = A(0,R_1,R_2)$ with $0 \leq R_1< R_2\leq+\infty$ and let $f : A \to \mathcal{X}$
be a left slice hyperholomorphic function. There exists $\{f_n\}_{n \in \mathbb{Z}} \subset  \mathcal{X}$
such that
\begin{equation}
f(q) = \sum_{n \in \mathbb{Z}} q^n f_n
\end{equation}
for all $q \in A$.
\end{Tm}
We now discuss the isolated real singularities. Let $U \subseteq \hh$ be an axially symmetric open set
 that intersects the real line and let $f:U\to \mathcal{X}$.
A point $\alpha\in  U \cap \mathbb{R}$ is called isolated singularity of $f$ if $\alpha \cup  U$ is open and if $\alpha \not\in U$.

Let $\alpha \in \mathbb{R}$ and  $U$ be a neighborhood of $\alpha$. Let $f:U\setminus \{\alpha\}\to \mathcal{X}$ be a slice left (or right) slice hyperholomorphic function and the open set $A(\alpha, 0, \varepsilon)$ defined by $0<|p-\alpha|<\varepsilon$ which is contained in $U$ for suitable $\varepsilon >0$.
By the previous result, there exists a uniquely determined series
$$
\sum_{n \in \mathbb{Z}} (q-\alpha)^n f_n
$$
which converges in for all $p$ such that  $0<|p-\alpha|<\varepsilon$
$$
f(q) = \sum_{n \in \mathbb{Z}} (q-\alpha)^n f_n.
$$
This is the Laurent expansion of $f$ at $\alpha\in \mathbb{R}$ and the element $f_{-1}$ is called the residue of $f$ at $\alpha$. The following definition fixes the terminology.
\begin{Dn} The isolated singularity $\alpha\in \mathbb{R}$ is called removable singularity of $f$ if $f_n=0$ for all negative integer $n$. If moreover, $f_0=0$, then $\alpha$ is a zero of $f$. If $\alpha$ is a zero of $f$ and $f$ is not identically zero in a neighborhood
of $\alpha$ then the smallest positive integer $n$ such that $f_{n}\not=0$ is called order of the zero $\alpha$.
The isolated singularity $\alpha$ is called a pole of $f$ if there exists a negative integer $m$ such that
 $f_{m}\not=0$ and $f_{n}\not=0$ for all integers $n\leq m-1$. The integer $m$ is called the order of the pole  $\alpha$.\\
 If $\alpha$ is not a removable singularity and is not a pole of $f$, then  $\alpha$ is called an essential singularity of $f$.
\end{Dn}
\begin{Tm}[Riemann's theorem of removability of singularities]\label{1.10.3}
Let $U\subseteq \hh$ be an axially symmetric open set that intersects the real line and let $\mathcal{X}$ be a Banach space.
Let $f: U \to \mathcal{X}$ be a left slice hyperholomorphic function and $\alpha \in \mathbb{R}$ be an isolated singularity of $f$ defined in a neighborhood of $\alpha$. If $f$ is bounded then $\alpha$ is a removable singularity, see \cite{GSS_book}.
\end{Tm}
\begin{proof}
It follows as in the case of scalar valued functions in the case $\alpha\in \mathbb{R}$.
\end{proof}
\begin{Tm}\label{1.10.4}
Let $U\subseteq \hh$ be an axially symmetric open set that intersects the real line and let $\mathcal{X}$ be a Banach space.
Let $f: U \to \mathcal{X}$ be a left slice hyperholomorphic function and $\alpha \in \mathbb{R}$ be an isolated singularity of $f$ defined in a neighborhood of $\alpha$. Then
the Laurent series of $f$ at $\alpha \in \mathbb{R}$ is of the form
$$
f(p)=\sum_{n=N}^\infty (q-\alpha)^n f_n, \ \ for \ \ f_N\not=0
$$
if and only if there exist constants $C<\infty$ and $c>0$ such that for some $\varepsilon>0$ we have
$$
c|p-\alpha|^N \leq \|f(p)\|\leq |p-\alpha|^N,\ \ \      0<|p-\alpha|<\varepsilon.
$$
\end{Tm}
\begin{proof}
It follows as in the case of scalar valued functions.
\end{proof}

\section{Additive splitting through Cauchy transform}\label{tre}
Since the Cauchy formulas of slice hyperholomorphic functions are based on an integral computed on complex planes $\mathbb{C}_j$, for $j\in \mathbb{S}$, the contours of integration have to be taken on open sets $U\subset \mathbb{H}$ intersected with the complex plane $\mathbb{C}_j$.
On a given complex plane $\mathbb{C}_j$, $j\in \mathbb{S}$, the definition of contour as well as the definition below are the same as in the complex case, but we repeat them for completeness.
\begin{Dn}[$\mathcal{C}^1$-contour on $\mathbb{C}_j$]\label{def31}
A set $\Gamma_j \subseteq \mathbb{C}_j$, $j\in \mathbb{S}$, is called a connected $\mathcal{C}^1$-contour on $\mathbb{C}_j$ if there exist real numbers $a,b$ with $a<b$ and a $\mathcal{C}^1$-function
$\gamma_j:[a,b]\to \mathbb{C}_j$ with $\Gamma_j=\gamma_j([a,b])$ such that
\begin{itemize}
\item[(i)] $\gamma_j'(t)\not=0$, for all $a\leq t\leq b$,
\item[(ii)] $\gamma_j(t)\not=\gamma_j(\tau)$ for all $a\leq t <\tau<b$,
\item[(iii)] it is $\gamma_j(b)\not=\gamma_j(t)$ for all $a\leq t<b$ or $\gamma_j(b)=\gamma_j(a)$ with $\gamma_j'(b)=\gamma_j'(a)$.
\end{itemize}
Moreover:
in the case the $\mathcal{C}^1$-parametrization $\gamma_j$ of $\Gamma_j$ is such that $\gamma_j(b)=\gamma_j(a)$ then $\Gamma_j$ is called closed contour.
\\
The union of a finite number of pairwise disjoint connected $\mathcal{C}^1$-contours is called $\mathcal{C}^1$-contour.
\end{Dn}
\begin{Dn}[Piecewise $\mathcal{C}^1$-contour on $\mathbb{C}_j$]
A set $\Gamma_j \subseteq \mathbb{C}_j$, for $j\in \mathbb{S}$, is called a connected piecewise $\mathcal{C}^1$-contour on $\mathbb{C}_j$ if one of the following conditions hold:
\begin{itemize}
\item[(i)]
There exist real numbers $a<b$ and a $\mathcal{C}^1$-function $\gamma_j:[a,b]\to \mathbb{C}_j$ with $\Gamma_j=\gamma_j([a,b])$ such that:
conditions (i) and (ii) in Definition \ref{def31} hold and when $\gamma_j(b)=\gamma_j(a)$ it is $\gamma_j'(a)/\gamma_j'(b)\in \mathbb{C}_j\setminus (-\infty,0]$.
\item[(ii)] For finitely many real numbers $t_\ell$, $\ell=1,...,m$ such that $a=t_1<...<t_m=b$ there exists a continuous function $\gamma_j:[a,b]\to \mathbb{C}_j$ with
$\Gamma_j=\gamma_j([a,b])$ such that:
\begin{itemize}
 \item[(a)]for each $\ell \in [1,m-1]$, the function $\gamma_{j,\ell}=\gamma_j\Big|_{[t_\ell,t_{\ell+1}]}$ belongs to $\mathcal{C}^1([t_\ell,t_{\ell+1}])$ and $\gamma_{j,\ell}'(t)\not=0$. Moreover, we assume that
      $\gamma_{j,\ell}(t_{\ell+1})/\gamma_{j,\ell+1}(t_{\ell+1})\in \mathbb{C}_j\setminus (-\infty,0]$ for each $\ell \in [1,m-2]$
     \item[(b)] $\gamma_j(t)\not=\gamma(\tau)$ for all $a\leq t< \tau\leq b$.
\end{itemize}
\item[(iii)]
Condition (a) holds and
\begin{itemize}
\item[(b')] $\gamma_j(t)\not=\gamma_j(\tau)$  for all $a\leq t< \tau< b$
and $\gamma_j(a)=\gamma_j(b)$ with $\gamma_j'(a)/\gamma_j'\in  \mathbb{C}_j\setminus (-\infty,0].$
\end{itemize}
\end{itemize}
\end{Dn}
A not necessarily connected piecewise $\mathcal{C}^1$-contour is the union of a finite number of pairwise disjoint connected piecewise $\mathcal{C}^1$-contours. Such a contour is closed if each connected component is closed. The orientation has to be considered as for contours in the complex plane.
A not necessarily connected piecewise $\mathbb{C}^1$-contour $\Gamma$ is called oriented if on each connected component of $\Gamma$ an orientation is fixed.

\begin{Tm}[Cauchy's integral theorem]
\label{CauchyThm}
Let $U\subset \H$ be an open set and let $\mathcal{X}$ be a two-sided quaternionic Banach space. Let
$f:\overline{U}\to \mathcal{X}$  be a continuous functions such that
$f\in\mathcal{SH}_{L}(U,\mathcal{X})$.
Moreover, let $ D_j\subset U\cap\mathbb{C}_j$ for all $j \in\mathbb{S}$  be an open and bounded subset of the complex plane $\mathbb{C}_j$ with $\overline{D_j}\subset U\cap\mathbb{C}_j$ such that  $\partial D_j$ is a piecewise $\mathcal{C}^1$-contour. Then
\begin{equation*}
\int_{\partial D_j} d s _j\, f( s ) = 0,
\end{equation*}
where $ds_{ j} = -j ds$.
\end{Tm}
As we mentioned already in the previous section, the Cauchy kernel to be used in the Cauchy formula for left slice hyperholomorphic functions is
$$
 S_L^{-1}(s,p):=-(p^2-2p{\rm Re}(s)+|s|^2)^{-1}(p-\overline{s}).
$$
It is a function slice hyperholomorphic on the left in the variable $p$ and on the right in the variable $s$.
In the case of right regular hyperholomorphic, the kernel  is
$$
S_R^{-1}(s,q):=-(q-\bar s)(q^2-2{\rm Re}(s)q+|s|^2)^{-1},
$$
which is slice hyperholomorphic on the right in the variable $q$ and on the left in $s$.
The Cauchy formula holds  for slice hyperholomorphic functions with values in a quaternionic Banach space:
\begin{Tm}[Cauchy formulas] \label{CauchygeneraleVV}
Let $\mathcal X$ be a  two sided quaternionic Banach algebra and let $W$  be an open set in $\hh$.
 Let $\overline{U}\subset W$ be an axially symmetric open set such that
 $\partial (U\cap \mathbb{C}_j)$ is a piecewise $\mathcal{C}^1$-contour  for every
$j\in\mathbb{S}$. Set  $ds_j=ds/ j$.
If  $f:W \to \mathcal{X}$ is a left slice hyperholomorphic, then, for $q\in U$, we have
\begin{equation}\label{cauchynuovo}
 f(p)=\frac{1}{2 \pi}\int_{\partial (U\cap \mathbb{C}_j)} S_L^{-1}(s,p)ds_j f(s),
\end{equation}
if  $f:W \to \mathcal{X}$ is a right slice hyperholomorphic, then, for $q\in U$, we have
\begin{equation}\label{Cauchyright}
 f(q)=\frac{1}{2 \pi}\int_{\partial (U\cap \mathbb{C}_j)}  f(s)ds_j {S}_R^{-1}(s,q),
\end{equation}
and the integrals (\ref{cauchynuovo}), (\ref{Cauchyright})  do not depend on the choice of the imaginary unit $j\in\mathbb{S}$ and on
$U\subset W$.
\end{Tm}


Before we state and prove our main results we need some more preliminaries.
\begin{Dn}
Let $V\in \mathbb{H}$ be an open set such that  $\mathbb{H}\setminus V$ is a bounded set.
We say that a function $f$ is slice $\mathcal{C}^k$-continuous, slice hyperholomorphic function
$f:V\bigcup \{\infty\}\to \mathcal{X}$ if $f$ is slice $\mathcal{C}^k$-continuous, slice hyperholomorphic function on $V$ and
$f(1/p)$ is slice $\mathcal{C}^k$-continuous, slice hyperholomorphic function on
$$
\{p\in \mathbb{H} \ :\ \frac{1}{p}\in V\}\cup\{0\}.
$$
\end{Dn}

\begin{Dn}[Splitting]
Let $U_+\subset\mathbb H$ be a bounded open set with piecewise
smooth boundary $\partial U_+ $ and let $U_-=\mathbb H\setminus \overline{U_+}$.
Let $f:\partial U \to \mathcal{X}$ be a slice continuous function.
We say that $f$ splits additively with respect to $\partial U_+$ if there exist functions
$f_{-}:\overline{U}_{-}\cup \{\infty\}\to \mathcal{X}$ and $f_{+}:\overline{U}_{+}\to \mathcal{X}$ where
$f_{-}$ is slice continuous on $\overline{U}_{-}$ and slice hyperholomorphic in ${U}_{-}\bigcup\{\infty\}$ and $f_+$ is slice continuous on $\overline{U}_+$ and slice hyperholomorphic in
$U_+$ such that, on $\partial U_+$ we have
$$
f=f_++f_-.
$$
\end{Dn}
We note that the splitting introduced in the previous definition is also called a global splitting.
We have the following results:
\begin{Tm}[Uniqueness of the splitting]
Let $\partial U_+$ be the smooth boundary of $U_+\subset \mathbb{H}$ and let $f:\partial U_+ \to \mathcal{X}$ be a slice continuous function.
Let $(f_-,f_+)$ be a global splitting of a continuous function $f$ with respect to $\partial U_+$. Then adding a constant we can always assume that
$f_-(\infty)=0$. With this condition the splitting $(f_-,f_+)$ is uniquely determined.
\end{Tm}
\begin{proof}
This assertion follows as in the complex case, see \cite{GL}, and it is a direct consequences of Liouville's theorem.
\end{proof}
We note that an example of splitting on the quaternionic unit sphere is given by the Wiener functions introduced in \cite{acks}, see also \cite{ACS_book}, Section 6.4.
\begin{Dn}\label{Cautrans}
Let $U_+\subset \mathbb H$ be a bounded, axially symmetric set such that
 $\partial (U_+\cap \mathbb{C}_j)$,  is the piecewise $\mathcal{C}^1$-contour  of $\partial(U_+\cap \mathbb{C}_j)$ for any $j\in \mathbb{S}$. Let  $f:\partial U_+ \to \mathcal{X}$ be a left slice continuous function.
We define, the left Cauchy integral transform with respect to  $\partial U_+$ of $f$ as
\begin{equation}\label{cauchyTL}
 \hat{f}(p)=\frac{1}{2 \pi}\int_{\partial (U_+\cap \mathbb{C}_j)} S_L^{-1}(s,p)ds_j f(s), \ \ \ p\in  \mathbb{H}\setminus \partial U_+.
\end{equation}
 We define, the right Cauchy integral transform with respect to  $\partial U_+$ of $f$ as
\begin{equation}\label{cauchyTR}
 \hat{f}(p)=\frac{1}{2 \pi}\int_{\partial (U_+\cap \mathbb{C}_j)}  f(s)ds_j {S}_R^{-1}(s,p), \ \ \ p\in  \mathbb{H}\setminus \partial U_+.
\end{equation}
\end{Dn}

\begin{La}\label{lemmaind} With the notation in Definition \ref{Cautrans}, the function defined in (\ref{cauchyTL}) does not depend on the choice of $j\in\mathbb S$ and
is left slice hyperholomorphic in $(\mathbb{H}\cup \{\infty\})\setminus \partial U_+$.
The function
defined in (\ref{cauchyTR}) does not depend on the choice of $j\in\mathbb S$ and
is right slice hyperholomorphic in $(\mathbb{H}\cup \{\infty\})\setminus \partial U_+$.
\end{La}
\begin{proof} We consider (\ref{cauchyTL}) since the argument for the function in \eqref{cauchyTR} is similar.
First of all, we note that the integral does not depend on the choice of the imaginary unit $j\in\mathbb S$, see Theorem 7.3.3 in \cite{ACS_book}. To prove that $f$ is slice hyperholomorphic we follow the reasoning in the proof of Theorem 7.1.3 in \cite{ACS_book} to show that $\hat f$ admits slice derivative for any $p\in \mathbb H\setminus \partial U$.
\end{proof}
An immediate consequence, proved with standard arguments is:
\begin{Cy}\label{derivative}
Under the assumptions of Lemma \ref{lemmaind} we have:
$$
 \hat{f}'(p)=\partial_{x_0} \hat{f}(p)
=\frac{1}{ \pi}
\int_{\partial (U\cap \mathbb{C}_j)}  (p^2-2s_0p+|s|^2)^{-2} (p-\overline{s})^{2*} ds_j f(s)
$$
where
\begin{equation}\label{stellina_qua}
(p-\overline{s})^{2*}=\sum_{k=0}^{2}\frac{2!}{(2-k)!k!} p^{2-k}\overline{s}^k,
\end{equation}
\end{Cy}

\begin{Tm}\label{3.2.due}
Let $U_+\subset\mathbb H$ be a bounded, axially symmetric open set. Assume that for any $j\in\mathbb S$ the set $\partial (U\cap \mathbb{C}_j)$ is piecewise $\mathcal{C}^1$. Let  $f:\partial U \to \mathcal{X}$ be a left slice continuous function.
Let us define
\begin{equation}\label{frfr}
 \hat{f}_+:=\hat{f}\Big|_{U_+},\ \ \ \hat{f}_-:=\hat{f}\Big|_{U_-\cup \{\infty\}},\ \ \
\end{equation}
be the two parts of the left Cauchy integral transform $\hat f$ of $f$. Then the two following conditions are equivalent:
\begin{itemize}
\item[(i)]
The function $f$ splits with respect to $\partial U$.
\item[(ii)]
The function $\hat{f}_+$ admits a continuous extension to $\overline{U}_+$ and  $\hat{f}_-$ admits a continuous extension to $\overline{U}_-$.
\end{itemize}
Moreover, this is the splitting of $f$ which vanish at infinity, namely $f=\hat{f}_++(-\hat{f}_-)$.
\end{Tm}
\begin{proof}
First of all we note that to say that $f:\partial U \to \mathcal{X}$ is a left slice continuous function trivially implies that $f:\partial U\cap\mathbb C_j \to \mathcal{X}$ is continuous for all $j\in\mathbb S$ and the integral defining $\hat f$ does not depend on $j\in\mathbb S$ by Lemma \ref{lemmaind}.

Step $(i)\Rightarrow (ii)$.
Let $f=f_++f_-$ be the splitting of $f$ with the condition $f_-(\infty)=0$.
Since $f_-$ is left slice hyperholomorphic at infinity and $f_-(\infty)=0$ we have thanks to Theorem \ref{1.10.4} that
$$
f_-(s)=O\Big(\frac{1}{|s|}\Big) \ \ \text{as}  \ \ |s|\to\infty.
$$
Using the Cauchy integral theorem and the estimate in Theorem \ref{1.3.6}
we have that
\begin{equation}
\label{323}
\lim_{r\to\infty}
 \int_{\partial (B_r(0)\cap \mathbb{C}_j)} S_L^{-1}(s,p)ds_j f_-(s)=0, \ \ \text{for every }\ \  p\in  \mathbb{H}\setminus \partial U
\end{equation}
where $B_r(0)$ is the ball in $\mathbb{H}$ centered at the origin and of radius $r>0$.

To prove that $f_+=\hat{f}_+$ in $U_+$ and $f_-=\hat{f}_-$ in $U_-$
we first let $p\in U_+$. If $r\in (0,\infty)$ is so large that $U_+$ in contained in the ball $B_r(0)$ by the Cauchy integral theorem we have
$$
 \int_{\partial (U_+\cap \mathbb{C}_j)} S_L^{-1}(s,p)ds_j f_-(s)= \int_{\partial (B_r(0)\cap \mathbb{C}_j)} S_L^{-1}(s,p)ds_j f_-(s).
$$
The relation (\ref{323}) shows that
$$
 \int_{\partial (U_+\cap \mathbb{C}_j)} S_L^{-1}(s,p)ds_j f_-(s)=0.
$$
Since  $f=f_++f_-$ on $\partial (U_+\cap \mathbb{C}_j)$ for  $j\in \mathbb{S}$ by the Cauchy formula we get
$$
 f_+(p)=\frac{1}{2 \pi}\int_{\partial (U\cap \mathbb{C}_j)} S_L^{-1}(s,p)ds_j f_+(s),
$$
and this implies
$$
 f_+(p)=\frac{1}{2 \pi}\int_{\partial (U\cap \mathbb{C}_j)} S_L^{-1}(s,p)ds_j f(s)=\hat{f}_+(p).
$$
We now take $p\in U_-$ and we observe that by the Cauchy theorem we have
\begin{equation}
\label{324}
\int_{\partial (U_-\cap \mathbb{C}_j)} S_L^{-1}(s,p)ds_j f_+(s)=0.
\end{equation}
We first assume that $p$ belongs to a bounded component of $U_-$ and denote by $\partial (U_-\cap \mathbb{C}_j)_0$ the part of the boundary
of this bounded component (with the orientation induced by $\partial (U\cap \mathbb C_j)$).
By the Cauchy theorem we have
\begin{equation}
\label{325}
\int_{\partial (U_-\cap \mathbb{C}_j)\setminus \partial (U_-\cap \mathbb{C}_j)_0} S_L^{-1}(s,p)ds_j f_-(s)=0
\end{equation}
so by the Cauchy formula
it follows that
\begin{equation}
\label{326}
\int_{ \partial (U_-\cap \mathbb{C}_j)_0} S_L^{-1}(s,p)ds_j f_-(s)=-f_-(p).
\end{equation}
From (\ref{324})(\ref{325})(\ref{326}) and the definition of $\hat{f}_-(p)$ we get
\[
\begin{split}
\hat{f}_-(p)&=\int_{ \partial (U_-\cap \mathbb{C}_j)} S_L^{-1}(s,p)ds_j f(s)
\\
&
=
\int_{ \partial (U_-\cap \mathbb{C}_j)} S_L^{-1}(s,p)ds_j [f_+(s)+f_-(s)]
\\
&=-f_-(p).
\end{split}
\]
We now consider the case when $p$ belongs to the unbounded component of $U_-\cap \mathbb{C}_j$ whose boundary is denoted by $\partial (U_-\cap \mathbb{C}_j)_\infty$.
By the Cauchy integral theorem we have
\begin{equation}
\label{327}
\int_{\partial (U_-\cap \mathbb{C}_j)\setminus \partial (U_-\cap \mathbb{C}_j)_\infty} S_L^{-1}(s,p)ds_j f(s)=0 .
\end{equation}
The Cauchy formula, for $r$ large enough, yields
$$
f_-(p)=-\int_{ \partial (U_-\cap \mathbb{C}_j)_\infty} S_L^{-1}(s,p)ds_j f_-(s)+
\int_{ \partial (B_r(0)\cap \mathbb{C}_j)_\infty} S_L^{-1}(s,p)ds_j f_-(s)
$$
so the relations (\ref{323}) and (\ref{327}) imply
$$
f_-(p)=-\int_{ \partial (U_-\cap \mathbb{C}_j)_\infty} S_L^{-1}(s,p)ds_j f_-(s)=
-\int_{ \partial (U_-\cap \mathbb{C}_j)} S_L^{-1}(s,p)ds_j f_-(s).
$$
Now using (\ref{324}) together with the definition of $\hat{f}_-$ we obtain
\[
\begin{split}
f_-(p)&=-\int_{ \partial (U_-\cap \mathbb{C}_j)_\infty} S_L^{-1}(s,p)ds_j [f_+(s)+f_-(s)]
\\
&
=
-\int_{ \partial (U_-\cap \mathbb{C}_j)} S_L^{-1}(s,p)ds_j f(s)
\\
&
=-\hat{f}_-(p).
\end{split}
\]
 Conversely, we prove $(ii)\Rightarrow (i)$. We have to show that if we fix $q\in \partial U_+$  and $\varepsilon >0$ we have
\begin{equation}
\label{ciccio}
|\hat{f}_+(q)-\hat{f}_-(q)-f(q)|\leq \varepsilon
\end{equation}
where $\hat{f}_+(q)-\hat{f}_-(q)$ denote the continuous extensions to $\partial U_+$.  We denote by
$\Gamma_j:=\partial (U_+\cap \mathbb{C}_j)$, for  $j\in \mathbb{S}$.
Assume that $w_0\in \Gamma_{j,0}$ where $\Gamma_{j,0}$ is the connected component of $\Gamma_{j}$ and let
$\gamma_j:[a,b]\to \Gamma_{j,0}$ be a $\mathcal{C}^1$-parametrization of $\Gamma_{j,0}=\partial (U\cap \mathbb{C}_j)_0$.
We recall that the functions  $f$, $f_+$ and $f_-$ are continuous on $U$ and in particular on $\Gamma_j=\partial (U\cap \mathbb{C}_j)$ so we can assume that
$w_0=\gamma_j(t_0)$, where $a<t_0<b$ is a point such that $\gamma_j'(t_0)\not=0$.
So we can find $\delta >0$, $c\in (0,1)$  and two sequences $w_n^+\in U_+\cap \mathbb{C}_j$ and $w_n^-\in U_-\cap \mathbb{C}_j$
such that
\begin{equation}
\label{3.2.11}
w_n^+\to w_0,\ \ \ \  w_n^-\to w_0,\ \ \text{as} \ \ n\to\infty
\end{equation}
that, for $n\in \mathbb{N}^*$ satisfy the conditions:
\begin{equation}
\label{3.2.9}
c|t-t_0|\leq |\gamma_j(t)-w_0|,\ \ \ \ |t-t_0|<\delta
\end{equation}
\begin{equation}
\label{3.2.10}
|f(\gamma_j(t))-f(w_0)|\leq \frac{c^2\varepsilon}{32 \max\{|\gamma_j'(t)| \ :|t-t_0|<\delta\  \text{and} \ j\in \mathbb{S}\}},\ \ \ \ |t-t_0|<\delta
\end{equation}
\begin{equation}
\label{3.2.12}
\begin{split}
&
\frac{1}{2}|\gamma_j(t)-w_0|+\frac{1}{4}  |w_n^+-w_n^-| \leq    |\gamma_j(t)-w_n^+|,\ \ \ \ |t-t_0|<\delta
\\
&
\frac{1}{2}|\gamma_j(t)-w_0|+\frac{1}{4}  |w_n^+-w_n^-| \leq |\gamma_j(t)-w_n^-|,\ \ \ \ |t-t_0|<\delta
 \end{split}
\end{equation}
thank to (\ref{3.2.9}) the estimates in (\ref{3.2.12}) become
\begin{equation}
\label{3.2.13}
\begin{split}
&
\frac{c}{4}(|t-t_0|+  |w_n^+-w_n^-|) \leq   |\gamma(t)-w_n^+|,\ \ \ \ |t-t_0|<\delta
\\
&
\frac{c}{4}(|t-t_0|+  |w_n^+-w_n^-|) \leq |\gamma(t)-w_n^-|,\ \ \ \ |t-t_0|<\delta.
 \end{split}
\end{equation}
Since we have defined $\Gamma_{j,0}=\partial (U\cap \mathbb{C}_j)_0$ we set
\begin{equation}\label{cau0L}
 \hat{f}_{+,\Gamma_0}(w^+_n)=\frac{1}{2 \pi}\int_{\partial (U\cap \mathbb{C}_j)_0} S_L^{-1}(s,w^+_n)ds_j f(s),
\end{equation}
and
\begin{equation}\label{cau0R}
 \hat{f}_{-,\Gamma_0}(w^-_n)=\frac{1}{2 \pi}\int_{\partial (U\cap \mathbb{C}_j)_0} {S}_L^{-1}(s,w^-_n) ds_j f(s)
\end{equation}
so we obtain
\[
\begin{split}
(\hat{f}_{+}(w^+_n)&-\hat{f}_{-}(w^-_n))-(\hat{f}_{+,\Gamma_0}(w^+_n)-\hat{f}_{-,\Gamma_0}(w^-_n))
\\
&
=\frac{1}{2 \pi}\int_{\partial (U\cap \mathbb{C}_j)_0\setminus \partial (U\cap \mathbb{C}_j)_0} [S_L^{-1}(s,w^+_n)-S_L^{-1}(s,w^-_n)]ds_j f(s),
\end{split}
\]
but since by (\ref{3.2.11}) we have  $w_n^+\to w_0$ and $w_n^-\to w_0$ as $n\to\infty$ so
$$
S_L^{-1}(s,w^+_n)-S_L^{-1}(s,w^-_n)\to 0 \ \ \ \text{as} \ \ \ n\to\infty.
$$
This shows that there exists $n_\varepsilon$ such that for $n\geq n_\varepsilon$ we have
$$
|(\hat{f}_{+}(w^+_n)-\hat{f}_{-}(w^-_n))-(\hat{f}_{+,\Gamma_0}(w^+_n)-\hat{f}_{-,\Gamma_0}(w^-_n))|<\varepsilon
$$
that is (\ref{ciccio}). Now assume that $-\Gamma_j:=-\partial (U_-\cap \mathbb{C}_j)$, for  $j\in \mathbb{S}$ is the boundary of the unbounded component of
$U_-\cap \mathbb{C}_j$, for  $j\in \mathbb{S}$, for $n$ sufficiently large, by the Cauchy formula and the Cauchy integral theorem, we have
$$
 f(w_0)=\frac{1}{2 \pi}\int_{\partial (U\cap \mathbb{C}_j)_0} S_L^{-1}(s,w^+_n)ds_j f(s)
\ \ \ {\rm and} \ \ \
\frac{1}{2 \pi}\int_{\partial (U\cap \mathbb{C}_j)_0} {S}_L^{-1}(s,w^-_n) ds_j f(s)=0.
$$
If
$-\Gamma_{j,0}:=-\partial (U_-\cap \mathbb{C}_j)_0$, for  $j\in \mathbb{S}$ is the boundary of one of bounded components of
$U_-\cap \mathbb{C}_j$, for  $j\in \mathbb{S}$, for $n$ sufficiently large, it is
$$
 \frac{1}{2 \pi}\int_{\partial (U\cap \mathbb{C}_j)_0} S_L^{-1}(s,w^+_n)ds_j f(w_0)=0
\ \ \ \text{and} \ \ \
\frac{1}{2 \pi}\int_{\partial (U\cap \mathbb{C}_j)_0} {S}_L^{-1}(s,w^-_n) ds_j f(w_0)=-f_-(w_0).
$$
So in both the cases, when we subtract the above relations, and for $n$ sufficiently large  we get
$$
 \frac{1}{2 \pi}\int_{\partial (U\cap \mathbb{C}_j)_0} \Big(S_L^{-1}(s,w^+_n)-{S}_L^{-1}(s,w^-_n)\Big)ds_j f(w_0)=f_-(w_0).
$$
So we obtain
$$
\hat{f}_{+,\Gamma_{j,0}}(w^+_n)-\hat{f}_{-,\Gamma_{j,0}}(w^-_n)-f(w_0)=\frac{1}{2 \pi}\int_{\partial (U\cap \mathbb{C}_j)_0} \Big(S_L^{-1}(s,w^+_n)-{S}_L^{-1}(s,w^-_n)\Big)ds_j
(f(s)-f(w_0)).
$$
Now we integrate on the path, where $\gamma_j(t)=\partial (U\cap \mathbb{C}_j)_0$, $t\in [a,b]$ and $\gamma_j(t_0)=w_0$ for $t_0\in (a,b)$, to get
\[
\begin{split}
\hat{f}_{+,\Gamma_{j,0}}(w^+_n)&-\hat{f}_{-,\Gamma_0}(w^-_n)-f(w_0)
\\
&
=\frac{1}{2 \pi}\int_a^b \Big(S_L^{-1}(\gamma_j(t),w^+_n)-{S}_L^{-1}(\gamma_j(t),w^-_n)\Big)(-j)(f(\gamma_j(t))-f(w_0))\gamma_j'(t)dt
\end{split}
\]
and we consider the splitting
$$
\hat{f}_{+,\Gamma_0}(w^+_n)
-\hat{f}_{-,\Gamma_0}(w^-_n)-f(w_0)=\mathcal{J}_1(n)+\mathcal{J}_2(n)
$$
where we have set
$$
\mathcal{J}_1(n):=\frac{1}{2 \pi}\int_{|t-t_0|\geq \delta} a^b \Big(S_L^{-1}(\gamma_j(t),w^+_n)-{S}_L^{-1}(\gamma_j(t),w^-_n)\Big)(-j)(f(\gamma_j(t))-f(w_0))\gamma_j'(t)dt,
$$
$$
\mathcal{J}_2(n):=\frac{1}{2 \pi}\int_{|t-t_0|<\delta} \Big(S_L^{-1}(\gamma_j(t),w^+_n)-{S}_L^{-1}(\gamma_j(t),w^-_n)\Big)(-j)(f(\gamma_j(t))-f(w_0))\gamma_j'(t)dt.
$$
From (\ref{3.2.11}) we have that
$$
\mathcal{J}_1(n)\to 0,\ \ \ \ \text{as} \ \ \ \ n\to\infty.
$$
We have to show that also
$$
\mathcal{J}_2(n)\to 0,\ \ \ \ \text{as} \ \ \ \ n\to\infty.
$$
So
$$
\|\mathcal{J}_2(n)\|=\frac{1}{2 \pi}\int_{|t-t_0|< \delta} \|(f(\gamma_j(t))-f(w_0))\|
\|S_L^{-1}(\gamma_j(t),w^+_n)-{S}_L^{-1}(\gamma_j(t),w^-_n)\|  \gamma_j'(t)dt
$$
and from (\ref{3.2.10})
set
$$
C_\gamma:=\max\{|\gamma_j'(t)| \ :|t-t_0|<\delta\  \text{and} \ j\in \mathbb{S}\}
$$
 we get
$$
\|\mathcal{J}_2(n)\|=\frac{1}{2 \pi}\int_{|t-t_0|< \delta}\frac{c^2\varepsilon}{32C_\gamma}
\
|S_L^{-1}(\gamma_j(t),w^+_n)-{S}_L^{-1}(\gamma_j(t),w^-_n)\|  |\gamma_j'(t)|dt
$$
and also
$$
\|\mathcal{J}_2(n)\|=\frac{1}{2 \pi}\frac{c^2\varepsilon}{32}\int_{|t-t_0|< \delta}
\|S_L^{-1}(\gamma_j(t),w^+_n)-{S}_L^{-1}(\gamma_j(t),w^-_n)\| dt
$$
and thanks to (\ref{3.2.13}) we finally have that there exists a constant $C>0$ such that
$$
\|\mathcal{J}_2(n)\|\leq C\varepsilon.
$$
\end{proof}

\section{The case of H\"older continuous functions}\label{quattro}
In the following we work in H\"older spaces. In the quaternionic case, to introduce a notion of distance we have to take into account not only points but also spheres. Thus we need the following definition.
\begin{Dn}
Let $\Gamma$ be a path in a complex plane $\mathbb{C}_j$, for $j\in \mathbb{S}$. We define
the distance between:
\begin{itemize}
\item[(i)] Two points $w$ and $s\in \mathbb{H}$ as
$
{\rm dist}(w,s):=|w-p|,
$
\item[(ii)]
Two spheres $[w]$ and $[s]$ as
$
{\rm dist}([w],[s]):=\inf \{|v-p|,\ v\in[w],\, p\in [s]\}.
$
\item[(iii)]
 The point  $w\in \mathbb{H}$ and a curve $\gamma$ as
$
{\rm dist}(w,\gamma):=\inf \{|w-s|,\ s\in\gamma\, \}.
$
\item[(iv)]
The sphere $[w]$ and the curve $\gamma$ as
$
{\rm dist}([w], \gamma):=\inf \{|v-p|,\ s\in\gamma,\, v\in [w]\}.
$
\item[(v)]
The  sphere $[p]$ and the point $s\in \mathbb{H}$  as
$
{\rm dist}([p],s):=\inf \{|q-s|,\ q\in[p]\, \}.
$
\end{itemize}
\end{Dn}

\begin{La}\label{Lemdist}
Let $p$, $s\in \mathbb{H}$ be such that $s\not\in[p]$ and  define the function
$$
\varphi_s(p):=(p^2-2s_0p+|s|^2)^{-2} (p-\overline{s})^{2*}.
$$
Then we have:
\begin{equation}\label{Stim}
|\varphi_s(p)|\leq ({\rm dist}([p],s))^{-2}
\end{equation}
\end{La}
\begin{proof}
Recall that using the $\star$-inverse it is
$$
(p-s)^{-\star}=(p^2-2s_0p+|s|^2)^{-1}(p-\overline{s})
$$
and we observe that $(p-s)^{-\star}$ can also be written as
$$
(p-s)^{-\star}=((p-\overline{s})^{-1}p(p-\overline{s})-s)^{-1}
$$
so
\[
\begin{split}
\varphi_s(p)&=(p^2-2s_0p+|s|^2)^{-1}(p^2-2s_0p+|s|^2)^{-1} (p-\overline{s})\star (p-\overline{s})
\\
&
=(p^2-2s_0p+|s|^2)^{-1}(p-s)^{-\star}\star (p-\overline{s})
\\
&
=(p^2-2s_0p+|s|^2)^{-1}(p-\overline{s})\star (p-s)^{-\star}
\\
&
=(p-s)^{-\star}\star (p-s)^{-\star}.
\end{split}
\]
Now we observe that $\tilde{p}:=(p-\overline{s})^{-1}p(p-\overline{s})\in [p]$ and so
$$
|(p-s)^{-\star}|=|(\tilde{p}-s)^{-1}|\leq\max_{\tilde{p}\in [p]}|(\tilde{p}-s)^{-1}|\leq ({\rm dist}([p],s))^{-1}
$$
so, we finally get
$$
|\varphi_s(p)|=|(p-s)^{-\star}\star (p-s)^{-\star}|\leq \Big(\max_{\tilde{p}\in [p]}|(\tilde{p}-s)^{-1}|\Big)^{2}= ({\rm dist}([p],s))^{-2}.
$$
\end{proof}
\begin{Rk} \label{rm_dist} {\rm If $p,s\in\C_j$, then $|\varphi_s(p)|=|p-s|^{-2}$, as $\varphi_s(p)$ then reduces to $(p-s)^{-2}$.}\end{Rk}
\begin{Rk}{\rm
As $\varphi_s(p)$ is a left slice function in $p$, the Representation Formula holds hence, given $s\in \hh$, suppose $s\in\C_j$ and take $p\in \hh$ such that $s\not\in [p]$. We have
$$
|\varphi_s(p)|=\dfrac{1}{2}|(1-j_pj)\varphi_s(p_j)+(1+j_pj)\varphi_s(\overline{p_j})|\leq ||p_j-s|^{-2}+|\overline{p_j}-s|^{-2}|\leq 2{\rm dist}([p],s)^{-2}.
$$
}
\end{Rk}
Let $U$ be an axially symmetric open set in $\mathbb H$.
The set of slice continuous functions $f:U\subseteq\hh\to \mathcal{X}$, where $\mathcal{X}$ is a quaternionic Banach space, will be endowed with the norm
$$
\|f\|_{U,0}=\sup_{s\in U}\|f(s)\|
$$
where $\|f(s)\|_\mathcal{X}$ denotes the norm in the Banach space $\mathcal{X}$.
\begin{Dn}
H\"older slice continuous functions of order $\alpha\in(0,1)$ on an axially symmetric open set $U\subseteq\mathbb H$
are defined as those  slice continuous functions $f:U\subseteq\hh\to \mathcal{X}$, with $f(x+jy)=\alpha(x,y)+j\beta(x,y)$, such that $\alpha$ and $\beta$ are H\"older continuous on $\mathcal{U}=\{(x,y)\in\mathbb R^2\ : \ x+jy\in U, j\in\mathbb S\}$; let $|\cdot|_{U,\alpha}$ be the usual H\"older seminorm on $\mathcal{U}$, then we endow the space $\mathcal C^\alpha (U, \mathcal X)$ of H\"older slice continuous functions of order $\alpha$ with the norm
$$
\|f\|_{U,\alpha}:=\|f\|_{U,0}+|\alpha|_{U,\alpha}+|\beta|_{U,\alpha}\;.
$$
\end{Dn}
These definitions are extended in an obvious way in the case of a closed set, in particular $\partial U$.
\begin{La}\label{Lemnorme}  Let  $U$ be a bounded set in $\mathbb H$ with piecewise smooth boundary $\Gamma$; we consider any $j\in\mathbb S$ and $f\in \mathcal{C}^\alpha(\Gamma_j, \mathcal{X})$, where $\Gamma_j=\Gamma\cap {\mathbb C}_j$ for $\alpha(0,1)$
and we assume that
$$
{\rm dist}([p],\Gamma_j)\geq 1.
$$
Then we have the following estimates:
\begin{equation}\label{casocontm1}
 \|\hat{f}'(p)\|
\leq \frac{|\partial (U\cap \mathbb{C}_j)|}{ \pi}
\dfrac{\|f\|_{\Gamma_j,0}}{ {\rm dist}([p], \Gamma_j)^2},
\end{equation}
and
\begin{equation}\label{casoHoldm2}
 \|\hat{f}'(p)\|
\leq 2\frac{|\partial (U\cap \mathbb{C}_j)|}{ \pi} \|f\|_{\Gamma_j,\alpha}{\rm dist}([p], \Gamma_j)^{\alpha-1}.
\end{equation}
\end{La}
\begin{proof}
By Corollary \ref{derivative} we have
$$
 \hat{f}'(p)
=\frac{1}{ \pi}
\int_{\partial (U\cap \mathbb{C}_j)}  (p^2-2s_0p+|s|^2)^{-2} (p-\overline{s})^{2*} ds_j f(s).
$$
By Theorem \ref{1.3.6} and Lemma \ref{Lemdist}, we get
$$
 \|\hat{f}'(p)\|
\leq \frac{|\partial (U\cap \mathbb{C}_j)|}{ \pi}
\dfrac{\|f\|_{\Gamma_j,0}}{ {\rm dist}([p], \Gamma_j)^2}
$$
which is estimate (\ref{casocontm1}).
To prove estimate (\ref{casoHoldm2}) we recall that, if $p\in \C_j$, we have
$$
\|\hat{f}'(p)\|\leq \frac{1}{ \pi}|\partial (U\cap \mathbb{C}_j)| \|f\|_{\Gamma_j,\alpha}{\rm dist}(p, \Gamma_j)^{\alpha-1}
$$
by the estimate of Remark \ref{rm_dist} and usual techniques, using the fact that ${\rm dist}([p],\Gamma_j)={\rm dist}(p,\Gamma_j)\geq 1$.

Now, as $\hat{f}'$ is a left-slice function, we use the Representation formula to obtain, as in Remark \ref{rm_dist}, that
$$\|\hat{f}'(p)\|\leq 2\dfrac{|\partial (U\cap\C_j)|}{\pi}\|f\|_{\Gamma_j,\alpha}{\rm dist}([p],\Gamma_j)^{\alpha-1}$$
\end{proof}
\begin{Tm}\label{Lem331}
 Let $\Gamma=\partial U$. For any $j\in\mathbb S$, let $f\in \mathcal{C}^\alpha(\Gamma_j, \mathcal{X})$ for $\alpha(0,1)$, where $\Gamma_j=\Gamma\cap{\mathbb C}_j$. Then there exists  a  constant $K>0$ such that the Cauchy transform satisfies
\begin{equation}\label{Holestim}
\|\hat{f}(p)\|\leq K \|f\|_{\Gamma_j,\alpha}{\rm dist}([p], \Gamma_j)^{\alpha-1}\ \ \text{for  all}\ \  p\in \mathbb{H}\setminus \Gamma_j\;.
\end{equation}
\end{Tm}
\begin{proof}
From the proof of the Lemma \ref{Lemnorme} we see that that for every axially symmetric neighborhood $A_j$ of $\Gamma:=\partial( U\cap \mathbb{C}_j)$, $j\in \mathbb{S}$
 there exists a positive constant $C_{A_j}<\infty$ such that
\begin{equation}\label{HolestimV}
\|\hat{f}(p)\|\leq C_{A_j}\|f\|_{\Gamma_j, \alpha}{\rm dist}([p], \Gamma_j)^{\alpha-1}\ \ \text{for \ all}\ \  p\in U\setminus A_j \ \ \text{with} \ \ {\rm dist}([p], \Gamma_j)\leq 1.
\end{equation}
Let $\varepsilon>0$ and consider the
axially symmetric neighborhood of a point $s_0\in \Gamma_j$ define as
$$
A_{j,\varepsilon}(s_0)=\{p\in \mathbb{H} \ :\ {\rm dist}([p], s_0)< \varepsilon \}.
$$
We now have to prove that for every $s_0\in \Gamma_j$ there exists $\varepsilon_0>0$ and a constant
$K_j>0$, independent of $f$, such that such that for all $p\in A_{\varepsilon_0}(s_0)\setminus\Gamma_j$ we have
(\ref{Holestim}) holds.
Let $s_0\in \Gamma_{0,j}$ be given, where $\Gamma_{0,j}$ is a connected component of $\Gamma_{j}$.
Let us choose a $\mathcal{C}^1$-parametrization $\gamma_j$ of $\Gamma_{0,j}$ such that
$\gamma_j:[-3,3]\to \Gamma_{0,j}$, $j\in \mathbb{S}$ with $\gamma_j(0)=s_0$.
 From the definition of the  $\mathcal{C}^1$-parametrization we have that there exist
constants
 $c_1<C_1<\infty$ such that
\begin{equation}\label{conditiparam}
c_1|t-\tau|\leq |\gamma_j(t)-\gamma_j(\tau)|\leq C_1|t-\tau|,\ \ \ \text{for\ all}\ \ -2\leq t,\tau\leq 2.
\end{equation}
Let us take $\varepsilon_0>0$ such that the following conditions hold
\begin{equation}\label{conditiepso}
\varepsilon_0<2c_1,\ \ \ \ A_{3\varepsilon_0}(s_0)\cap(\Gamma_j\backslash\gamma_j([-1,1])=\emptyset, \ \ \text{for }\ \  j\in \mathbb{S}.
\end{equation}
Let us take $p\in A_{\varepsilon_0}(s_0)$ and set
$$
\varepsilon={\rm dist}([p], \Gamma_j)\ \  \text{for }\ \  j\in \mathbb{S}
$$
since $s_0\in \Gamma$, then $\varepsilon\leq\varepsilon_0$ take a point $s'\in\Gamma$ with
$$
|[p]-s'|={\rm dist}([p], \Gamma_j)=\varepsilon.
$$
Since $p\in  A_{\varepsilon_0}(s_0)$ and $\varepsilon\leq\varepsilon_0$ then $s'\in A_{2\varepsilon_0}(s_0)$.
Since $s'\in \Gamma_j \cap A_{2\varepsilon_0}(s_0)$, it follows from the condition
$\varepsilon_0<2c_1$ in (\ref{conditiepso}) that $s'\in \gamma_j([-1,1])$.
Now let $t'\in [-1,1]$ be the parameter with $\gamma_j(t')=s'$, this parameter is uniquely determined by the condition (\ref{conditiparam}). Recall that since $\Gamma_j$, $j\in \mathbb{S}$ are closed curves, by slice hyperholomorphicity it follows that
\begin{equation}
\frac{1}{ 2\pi}
\int_{\partial (U\cap \mathbb{C}_j)}  (p^2-2s_0p+|s|^2)^{-2} (p-\overline{s})^{2*} ds_jf(s') =0
\end{equation}
so we have
\begin{equation}\label{modiffsec}
\hat{f}'(p)=\frac{1}{ 2\pi}
\int_{\partial (U\cap \mathbb{C}_j)}  (p^2-2s_0p+|s|^2)^{-2} (p-\overline{s})^{2*} ds_j[f(s)-f(s')].
\end{equation}
Since $t'\in [-1,1]$ and thanks to the condition
$\varepsilon_0<2c_1$ in (\ref{conditiepso}) we have the inequalities
$$
\frac{2\varepsilon }{c_1}\leq \frac{2\varepsilon_0 }{c_1}<1
$$
from which we get
$$
-2\leq t'-\frac{2\varepsilon }{c_1}< t'+\frac{2\varepsilon }{c_1}<2.
$$
So we can split the integral (\ref{modiffsec}) as the sum of three terms:
$$
\hat{f}'(p)=J_{1}(p)+J_{2}(p)+J_{3}(p),
$$
where
\[
\begin{split}
J_{1}(p)&:=\frac{1}{ 2\pi}
\int_{\Gamma_\setminus\gamma_j([-2,2])}  \varphi_{\gamma_j(t)}(p) ds_j[f(s)-f(s')],
\\
J_{2}(p)&:=\frac{1}{ 2\pi}
\int_{ t'-\frac{2\varepsilon }{c_1}}^{ t'+\frac{2\varepsilon }{c_1}}  \varphi_{\gamma_j(t)}(p) \, \gamma_j'(t)(-j)[f(\gamma_j(t))-f(s')]dt,
\\
J_{3}(p):&=
\frac{1}{ 2\pi}
\int_{-2}^{ t'-\frac{2\varepsilon }{c_1}} \varphi_{\gamma_j(t)}(p) \, \gamma_j'(t)(-j)[f(\gamma_j(t))-f(s')]dt
\\
&
+
\frac{1}{ 2\pi}
\int_{t'+\frac{2\varepsilon }{c_1}}^2  \varphi_{\gamma_j(t)}(p) \, \gamma_j'(t)(-j)[f(\gamma_j(t))-f(s')]dt.
\end{split}
\]
where we have set
$$
\varphi_s(p):=(p^2-2s_0p+|s|^2)^{-2} (p-\overline{s})^{2*}.
$$
Since
$$
|[p]-s'|={\rm dist}([p], \Gamma_j)=\varepsilon\leq\varepsilon_0
$$
by $A_{3\varepsilon_0}(s_0)\cap(\Gamma_j\backslash\gamma_j([-1,1])=\emptyset$ in (\ref{conditiepso})
we have
$$
|s-s_0|\geq 3\varepsilon_0,\ \ \ \text{for}\ \ \ s\in \Gamma_j\setminus \gamma_j([-2,2])
$$
so we get the estimate for $J_1$ by Proposition \ref{1.3.6}
\begin{equation}\label{STIMAJ1}
\|J_1(p)\|\leq \frac{|\Gamma_j|}{4\pi \varepsilon_0^2}\max_{s\in \Gamma_j}\|f(s)\|
\leq \frac{|\Gamma_j|}{4\pi \varepsilon_0^2}\|f\|_{\Gamma,\alpha}
\leq \frac{|\Gamma_j|}{4\pi \varepsilon_0^2}\|f\|_{\Gamma,\alpha}\varepsilon^{\alpha-1}.
\end{equation}
for $\varepsilon\in (0,1)$.
Set
$$
C_{2,j}:=\max_{t\in [-2,2]}|\gamma_j'(t)|.
$$
Since $|\gamma_j(t)-[p]|\geq {\rm dist}([p],\Gamma_j)=\varepsilon$ for $t\in [-3,3]$ it follows that
$$
\|J_{2}(p)\|\leq \frac{C_2}{ 2\pi\varepsilon^2}
\int_{ t'-\frac{2\varepsilon }{c_1}}^{ t'+\frac{2\varepsilon }{c_1}}  \,
\|f(\gamma_j(t))-f(s')\|dt
$$
since also on a slice
$$
\|f(\gamma_j(t))-f(s')\|\leq \|f\|_{\Gamma_j,\alpha}|\gamma_j(t)-s'|^\alpha,\ \ \ t\in [-3,3]
$$
we have
$$
\|J_{2}(p)\|\leq \frac{C_2C_1\|f\|_{\Gamma_j,\alpha}}{ 2\pi\varepsilon^2}
\int_{ t'-\frac{2\varepsilon }{c_1}}^{ t'+\frac{2\varepsilon }{c_1}}  \,
|t-\tau|^\alpha dt
$$
with some computations
$$
\|J_{2}(p)\|\leq \frac{2C_1C_2}{ c_1^2}
\|f\|_{\Gamma_j,\alpha}\varepsilon^{\alpha-1}.
$$
Finally, let us turn our attention to $J_3$. For
$$
|t-\tau|\geq \frac{2\varepsilon}{c_1}
$$
from (\ref{conditiparam}) we have
$$
\frac{1}{2}|\gamma_j(t)-\gamma_j(\tau)|\geq \frac{1}{2}|t-\tau| \geq \varepsilon.
$$
Since
$$
|\gamma_j(t)-[p]|=|s'-[p]|=\varepsilon
$$
and for
$$
|t-\tau|\geq \frac{2\varepsilon}{c_1}
$$
this implies
$$
|\gamma_j(t)-[p]|\geq |\gamma_j(t)-\gamma_j(t')|-|\gamma_j(t')-[p]|\geq |\gamma_j(t)-\gamma_j(t')|-\varepsilon
\geq \frac{1}{2}|\gamma_j(t)-\gamma_j(t')|
$$
so we get
\[
\begin{split}
\|J_{3}(p)\|&\leq
\int_{-2}^{ t'-\frac{2\varepsilon }{c_1}}
|\gamma_j'(t)|\, \frac{|f(\gamma_j(t))-f(s')|}{\|\gamma_j(t)-\gamma_j(t')\|^2}dt
\\
&
+
\int_{t'+\frac{2\varepsilon }{c_1}}^2  |\gamma_j'(t)|\, \frac{\|f(\gamma_j(t))-f(s')\|}{|\gamma_j(t)-\gamma_j(t')|^2}dt
\end{split}
\]
using now the H\"older continuity
$$
\|f(\gamma_j(t))-f(s')\|=
\|f(\gamma_j(t))-f(\gamma_j(t'))\|\leq \|f\|_{\Gamma_j,\alpha}|\gamma_j(t)-\gamma_j(t')|^\alpha
$$
with some computations recalling the definition of $C_2$ we get the estimate
$$
\|J_{3}(p)\|\leq \frac{2C_1C_2 c_1^{1-\alpha}}{1-\alpha}
\|f\|_{\Gamma_j,\alpha}\varepsilon^{\alpha-1}
$$
so taking
$$
K_j:=\frac{|\Gamma_j|}{4\pi \varepsilon_0^2}+\frac{2C_1C_2}{ c_1^2}+\frac{2C_1C_2 c_1^{1-\alpha}}{1-\alpha}
$$
we observe that the constant $K_j$ depends continuously on $j\in \mathbb{S}$, but since $\mathbb{S}$ is compact
we get the statement.
\end{proof}
We are now in the position to state and prove the main result.

\begin{Tm}
Let $U\subset \mathbb{H}$ be a bounded axially symmetric set and let $\Gamma:=\partial U_+$ be the piecewise $\mathcal{C}^1$-contour  of $U_+$.
We consider any $j\in\mathbb S$ and  $f\in \mathcal{C}^\alpha(\Gamma_j, \mathcal{X})$ for $\alpha\in (0,1)$, where $\Gamma_j=\Gamma\cap{\mathbb C}_j$.
Set:
\begin{equation}\label{frpp}
 \hat{f}_+:=\hat{f}\Big|_{U_+},\ \ \ \hat{f}_-:=-\hat{f}\Big|_{U_-},
\end{equation}
i.e. the two parts of the left Cauchy integral transform of $f$. Then the two following conditions are equivalent:
\begin{itemize}
\item[(i)]
The function $\hat{f}_+$ admits an H\"older continuous extension of order $\alpha\in (0,1)$ to $\overline{U}_+$,
and the function $\hat{f}_-$ admits an H\"older continuous extension of order $\alpha\in (0,1)$ to $\overline{U}_+\cup\{\infty\}$.
\item[(ii)]
 If we denote these extensions with the notation $f_+$ and $f_-$, then we have $f=f_++f_-$ on $\Gamma_j$.
\end{itemize}
\end{Tm}
\begin{proof}
Point (ii) follows from Theorem \ref{3.2.due}.
To prove  point (i) it is sufficient to show that there exists an axially symmetric neighborhood $A_{\Gamma_j}$ of the curve $\Gamma_j=\Gamma\cap\mathbb C_j$ and a constant $C_j>0$
such that
\begin{equation}\label{stimaGjpiu}
\|f_+(w_1)-f_+(w_2)\|\leq C_j\|f\|_{\Gamma_j}|w_1-w_2|
\end{equation}
for all $w_1,w_2\in A_{\Gamma_j}\cap (U_+\cap\mathbb{C}_j)$
and
\begin{equation}\label{stimaGjmeno}
\|f_-(w_1)-f_-(w_2)\|\leq C_j\|f\|_{\Gamma_j}|w_1-w_2|
\end{equation}
for all $A_{\Gamma_j}\cap (U_-\cap\mathbb{C}_j)$.
Since the proofs of the above estimates are the same we will consider just (\ref{stimaGjpiu}).
 Let us take an arbitrary  point $s_0\in \Gamma_j$ define as
$$
A_{j,\varepsilon}(s_0)=\{p\in \mathbb{H} \ :\ {\rm dist}([p], s_0)< \varepsilon \}.
$$
To prove our statement we have to find constants $\varepsilon_0>0$ and $C_j>0$ such that
(\ref{stimaGjpiu}) holds true for all $w_1,w_2\in  A_{j,\varepsilon}(s_0)\cap \mathbb{C}_j$, $j\in \mathbb{S}$.
By Theorem \ref{Lem331} there exists a constant $C_{0,j}>0$ such that
\begin{equation}\label{Holestimpp}
\|\hat{f}_+'(p)\|\leq C_{0,j} \|f\|_{\Gamma_j,\alpha}{\rm dist}([p], \Gamma_j)^{\alpha-1}\ \ \text{for  all}\ \  p\in U_+,
\end{equation}
Since $\Gamma_j$ is  $\mathcal{C}^1$-piecewise it is possible to find $\varepsilon$, $c\in (0,1)$ and a quaternion $v\in \mathbb{C}_j$ with $|v|=1$ such that $w\in A_{j,\varepsilon_0}(s_0)\bigcap U_+$ it is
\begin{equation}\label{3315}
A_{j,ct}(w+tv) \subseteq U_+,\ \ \ \text{for}\ \ t\in [0,\varepsilon_1].
\end{equation}
Let us set
\begin{equation}\label{seet}
\varepsilon_0=c\varepsilon_1 /4  \ \ \ \text{and} \ \ \ C_j=C_{0,j}\Big(1+\frac{2^{\alpha+1}}{c\alpha}\Big).
\end{equation}
We now show that with the constants chosen as in  (\ref{seet})  estimate (\ref{stimaGjpiu})
holds for all $w_1,w_2\in A_{j,\varepsilon_0}(s_0)\cap (U_+\cap\mathbb{C}_j)$.
Let $w_1,w_2\in A_{j,\varepsilon_0}(s_0)\cap (U_+\cap\mathbb{C}_j)$ be given.
We define
$\varepsilon=|w_1-w_2|$. Then
$$
2\varepsilon/c\leq 4\varepsilon_0/c=\varepsilon_1
$$
and it follows from (\ref{seet}) that
$$
w_1+tv\in U_+\ \ \ \  \text{and} \ \ \ w_2+tv\in U_+,\ \ \text{for all}\ t\in [0,2\varepsilon/c]
$$
and
$$
(1-t)(w_1+2v\varepsilon/c)+ t(w_2+2v\varepsilon/c)\in U_+,\ \ \text{for all}\ t\in [0,2\varepsilon/c],
$$
so we get
\[
\begin{split}
\|f_+(w_1)-f_+(w_2)\|&\leq \| f_+(w_1)-f_+(w_1+2v\varepsilon/c)\|+\| f_+(w_2)-f_+(w_2+2v\varepsilon/c)\|
\\
&
+\| f_+(w_1+2v\varepsilon/c)-f_+(w_2+2v\varepsilon/c)\|
\\
&
=\|\int_0^{2\varepsilon/c}f_+'(w_1+tv)vdt\|+\|\int_0^{2\varepsilon/c}f_+'(w_2+tv)vdt\|
\\
&
+
\|\int_0^{2\varepsilon/c}f_+'((1-t)(w_1+2v\varepsilon/c)+ t(w_2+2v\varepsilon/c))(w_2-w_1)dt\|
\end{split}
\]
but we have set $\varepsilon=|w_1-w_2|$ and $|v|=1$ this implies
\begin{equation}
\label{3316}
\begin{split}
\|f_+(w_1)-f_+(w_2)\|&\leq \|\int_0^{2\varepsilon/c}f_+'(w_1+tv)vdt\|+\|\int_0^{2\varepsilon/c}f_+'(w_2+tv)vdt\|
\\
&
+
\varepsilon \int_0^{2\varepsilon/c} \|f_+'((1-t)(w_1+2v\varepsilon/c)+ t(w_2+2v\varepsilon/c))\|dt.
\end{split}
\end{equation}
From (\ref{3315}) and (\ref{Holestimpp}) we have that
$$
\|\hat{f}_+'(w_\ell+tv)\|\leq C_{0,j} \|f\|_{\Gamma_j,\alpha}(ct)^{\alpha-1}\ \ \text{for  all}\ \  t\in [0,2\varepsilon /c],\ \ \ \ell=1,2,
$$
so we obtain
\begin{equation}\label{3317}
\int_0^{2\varepsilon /c}\|\hat{f}_+'(w_\ell+tv)\|dt  \leq C_{0,j} \|f\|_{\Gamma_j,\alpha} \int_0^{2\varepsilon /c}(ct)^{\alpha-1}dt=
\frac{2^\alpha C_{0,j}}{c\alpha} \|f\|_{\Gamma_j,\alpha}\varepsilon^\alpha  ,\ \ \ \ell=1,2.
\end{equation}
Observe that from (\ref{3315}) we also have
$$
A_{j,2\varepsilon}(w_2+2v\varepsilon/c)\subseteq U_+.
$$
Since
$$
\|(w_1+2v\varepsilon/c)-(w_2+2v\varepsilon/c)\|=|w_1-w_2|=\varepsilon
$$
it follows
$$
(1-t)(w_1+2v\varepsilon/c)+t(w_2+2v\varepsilon/c)\in \overline{A_{j,2\varepsilon}}(w_2+2v\varepsilon/c), \ \ \ t\in [0,1]
$$
this gives
$$
(1-t)(w_1+2v\varepsilon/c)+t(w_2+2v\varepsilon/c)\in U_+, \ \ \ t\in [0,1]
$$
and
$$
{\rm dist}\Big((1-t)(w_1+2v\varepsilon/c)+t(w_2+2v\varepsilon/c), \Gamma_j \Big)\geq \varepsilon
$$
so with (\ref{Holestimpp}) we obtain
\begin{equation}\label{3318}
\varepsilon \int_0^{2\varepsilon/c} \|f_+'((1-t)(w_1+2v\varepsilon/c)+ t(w_2+2v\varepsilon/c))\|dt\leq C_{0,j}
 \|f\|_{\Gamma_j,\alpha}\varepsilon^\alpha,
\end{equation}
we conclude saying that estimate (\ref{stimaGjpiu}) follows form estimates (\ref{3316}) (\ref{3317}) and (\ref{3318})
\end{proof}

\section{Fundamental solution of the global operator of slice hyperholomorphic functions}\label{cinque}

In the quaternionic setting, the Cauchy-Riemann operator is replaced by the operator $G_L$ (resp. $G_R$), introduced in \cite{GLOBAL}; when restricted to a slice $\C_j$, it reduces to the operator $y^2\overline{\partial}_j$.

 The left quaternionic global Cauchy-Riemann operator is given by
 $$G_Lf(q)=|\underline{q}|^2\dfrac{\partial f}{\partial x_0}(q)+\underline{q}\sum_{j=1}^3 x_j\dfrac{\partial f}{\partial x_j}(q)$$
 and, similarly, the right operator is given by
$$G_Rf(q)=|\underline{q}|^2\dfrac{\partial f}{\partial x_0}(q)+\sum_{j=1}^3 x_j\dfrac{\partial f}{\partial x_j}(q)\underline{q}\;.$$

\begin{Rk}\label{rmk1}
{\rm If $U\subset\hh$ is an axially symmetric open set and $f:U\to\hh$ is a continuous function such that $G_L(f)=0$ in the distribution sense, then
$$
v^2\left(\dfrac{\partial}{\partial u}+j\dfrac{\partial}{\partial v}\right)f(u+jv)=0
$$
for every $j\in\mathbb{S}$. Away from the real axis, this implies that $f(u+jv)$
is a holomorphic function; as it is also continuous, by a standard argument of one complex variable, we deduce that it is holomorphic on $U\cap\C_j$, hence by the arbitrarity of $j\in\mathbb S$, $f$ is slice hyperholomorphic.}
\end{Rk}

Given an axially symmetric open set  $U$ and a left slice function $V:U\to\hh$, we consider the equation
\begin{equation}\label{HCR}G_L(f)(q)=|\underline{q}|^2V(q)\;.\end{equation}
In order to study its solvability, we introduce an appropriate class of distributions.

\begin{Dn}Let $U\subseteq\hh$ be an axially symmetric open set
 and let $\mathcal{U}\subseteq\mathbb{R}\times \mathbb{R}$ be such that $p=u+j v\in U$ for all $(u,v)\in\mathcal{U}$ and all $j\in\mathbb{S}$. A (left) slice-test function is a function of the form
 $$\varphi(p)=\varphi(u+j v)=f_0(u,v)+jf_1(u,v)$$
 with $f_0,f_1\in\mathcal{C}^\infty_c(\mathcal{U})$, satisfying $f_0(u,-v)=f_0(u,v)$ and $f_1(u,-v)=-f_1(u,v)$. A right slice-test function is defined similarly, with $\varphi(p)=f_0(u,v)+f_1(u,v)j$.

 We denote the space of left (resp. right) slice-test functions by $\mathcal{SD}_L(U)$ (resp. $\mathcal{SD}_R(U)$).\end{Dn}

 The space $\mathcal{SD}(U)$ has the usual Fréchet topology; a right (resp. left) slice-distribution is a linear map $L:\mathcal{SD}_L(U)\to\H$ (resp. $L:\mathcal{SD}_R(U)\to\H$) which is right-linear (resp. left-linear) and continuous, i.e. for every $K\Subset U$ there exist an integer $m\geq 0$ and a constant $C>0$ such that
 $$
 |L(\varphi)|\leq C\sup_{|\alpha|\leq m}
 \left(\sup_{p=u+j v\in U}\left|\left(\frac{\partial}{\partial u}\right)^{\alpha_1}\left(\frac{\partial}{\partial v}\right)^{\alpha_2}\varphi(p)\right|\right)
 $$
 for every $\varphi\in\mathcal{SD}_L(U)$ (resp. $\varphi\in\mathcal{SD}_R(U)$) with $\mathrm{supp}\varphi\subseteq K$, where $|\alpha|=\alpha_1+\alpha_2$.

 Given $\mathcal{U}\subseteq\R^2$, invariant under the isometry $(u,v)\mapsto (u,-v)$, and its axially symmetric completion $U=\{x+jy\ :\ (x,y)\in\mathcal U,\, j\in\mathbb S\}\subseteq\hh$, we define the operators $P_\pm:\mathcal{C}^\infty_c(\mathcal{U})\to\mathcal{C}^\infty_c(\mathcal{U})$ as
 $$P_\pm f(u,v)=\dfrac{f(u,v)\pm f(u,-v)}{2}\;.$$
 It is easy to see that both are linear and bounded with respect to the Fréchet topology. Given $f\in\mathcal{C}^\infty_c(\mathcal{U})$, we construct $T(f)\in\mathcal{SD}_L(U)$ by
 $$T(f)(q)=P_+f(q_0,|\underline{q}|)+j_qP_-f(q_0,|\underline{q}|)$$
 where $j_q=\underline{q}/|\underline{q}|$. Then, $T:\mathcal{C}^\infty_c(\mathcal{U})\to\mathcal{SD}_L(U)$ is easily seen to be linear and bounded with respect to the Fréchet topology.

On the other hand, a test function $\varphi\in\mathcal{SD}_L(U)$ corresponds to a pair $(f_0,f_1)$ of functions in $\mathcal{C}^\infty_c(\mathcal{U})$
 such that $f_0$ is even in $v$ and $f_1$ is odd in $v$; clearly, $T(f_0+f_1)=\phi$. Hence, $T$ is an isomorphism.
 Therefore, a right slice-distribution corresponds to a classical distribution on $\mathcal{U}$, under the identification given by $T$.

 We also note that the operators $P_{\pm}$ extend as continuous linear isomorphisms to $L^p$ spaces as well.

 \begin{Pn}\label{testprodotto}Let $\psi:U\to\H$ be a right slice-$L^1_\mathrm{loc}$ function; suppose that, for some $j\in\mathbb{S}$, we have
 $$\int_{U\cap\C_j}\psi(u+jv)\varphi(u+jv)dudv=0\qquad \forall \varphi\in\mathcal{SD}_L(U)\;.$$
 Then $\psi\equiv 0$ almost everywhere on $U$.\end{Pn}
 \begin{proof} Let us write $\varphi=T(f)$ and $\psi=T(g)$, then
 $$\int_{U\cap \C_j}\psi(u+jv)\varphi(u+jv)dudv
 $$
 $$=\int_{U\cap\C_j}(P_+g(u,v)P_+f(u,v)+P_+g(u,v)jP_-f(u,v)+$$
$$ \hspace{2cm} +P_-g(u,v)jP_+f(u,v)-P_-g(u,v)P_-f(u,v))dudv\;.$$
 Now, we note that $U\cap\C_j$ is symmetric with respect to the $x$ axis, whereas the functions
 $P_+g(u,v)jP_-f(u,v)$ and $P_-g(u,v)jP_+f(u,v)$
 are odd with respect to the $v$ variable, hence their integral on $U$ vanishes. Therefore, our integral reduces to
 $$\int_{U\cap\C_j}(P_+g(u,v)P_+f(u,v)-P_-g(u,v)P_-f(u,v))dudv,$$
 which has to vanish for every choice of $f$ smooth and compactly supported. We notice that also the function
 $$P_+gP_-f - P_-gP_+f$$
 is odd with respect to the $v$ variable, so its integral on $U\cap\C_j$ vanishes, so
 $$\int_{U\cap\C_j}(P_+g(u,v)P_+f(u,v)-P_-g(u,v)P_-f(u,v))dudv$$
 $$=\int_{U\cap\C_j}(P_+g(u,v)P_+f(u,v)+P_+g(u,v)P_-f(u,v) -$$
 $$\hspace{2cm}- P_-g(u,v)P_+f(u,v)-P_-g(u,v)P_-f(u,v))dudv$$
 $$=\int_{U\cap\C_j}(P_+g(u,v)-P_-g(u,v))f(u,v)dudv,$$
 where we used that $P_++P_-=Id$. This integral must vanish for every $f$ smooth and compactly supported, so $P_+g-P_-g=0$ a.e., so $P_+g=P_-g$, but these two functions are respectively even and odd, then $P_+g=P_-g=0$ a.e..
 This easily implies that $g$ must vanish almost everywhere.\end{proof}

The previous result implies that two slice-$L^1_\loc$ functions give (by integration) the same slice-distribution if and only if they coincide almost everywhere.

\begin{Rk} {\rm Obviously, two slice-$L^1_\loc$ functions are also $L^1_\loc$ in the classical sense on $U$, so they give the same (classical) distribution on $U$, integrating with respect to the Lebesgue measure, if and only if they coincide almost everywhere. Indeed, the computations of the previous proposition also show that, in case of slice-functions, their action as (classical) distributions can be recovered from their action on a slice.

More explicitly, we consider the map
$$\pi:\mathcal{U}\times\mathbb{S}\to U$$
given by $\pi(u,v,j)=u+jv$; if we consider on $U$ the standard volume form $\mu=dx_0\wedge dx_1\wedge dx_2\wedge dx_3$ given by the inclusion of $U$ in $\R^4$, then
$$\pi^*\mu=v^2du\wedge dv\wedge d\sigma(j)$$
where $d\sigma$ is the standard volume form on the $2$-sphere.

Then, as a side-product of the computation performed in the previous proof and with the same notation, we obtain that
 \[
 \begin{split}
 \int_U\psi(q)f(q)dx_0dx_1dx_2dx_3&=\int_{\mathcal{U}\times\mathbb{S}}\!\!\!(P_+g(u,v)+P_-g(u,v)j)(P_+f(u,v)+jP_-f(u,v))v^2dudvd\sigma(j)
 \\
 &
 =\int_{\mathbb{S}}d\sigma\int_{\mathcal{U}}(P_+g(u,v)P_+f(u,v)-P_-g(u,v)P_-f(u,v))v^2dudv
 \\
 &
 =4\pi\int_{\mathcal{U}}(P_+g(u,v)P_+f(u,v)-P_-g(u,v)P_-f(u,v))v^2dudv.
 \end{split}
 \]
 The last integral is zero if and only if the function $P_+g(u,v)P_+f(u,v)-P_-g(u,v)P_-f(u,v)$ is zero almost everywhere. Therefore, we have
 $$\int_{U}\psi\varphi =0$$
 if and only if
 $$\int_{\mathcal{U}}(P_+\psi(u,v)P_+\varphi(u,v)-P_-\psi(u,v)P_-\varphi(u,v))dudv=0$$
 but this happens if and only if
 $$\int_{U\cap\C_j}\psi\varphi=0$$
 for some (and hence every) $j\in\mathbb{S}$.}\end{Rk}

 So, the choice of working on a slice or on the whole open set does not affect the result, as long as we are defining the space of distributions; however, such a choice becomes important when we want to define the adjoint of an operator. Namely, the \emph{global} (distributional) adjoint of $G_L$ is an operator $G_L^*$ which satisfies
 $$\int_{U}G^*_L(\varphi)(q)\psi(q)dq=\int_{U}\varphi(q)G_L(\psi)(q)dq$$
 for every $\varphi\in\mathcal{SD}_R(U)$ and $\psi\in\mathcal{SD}_L(U)$, where $dq$ is the standard $4$-dimensional volume measure on $U$. On the other hand, the \emph{slice} (distributional) adjoint of $G_L$ is an operator $G_L^{*s}$ which satisfies
$$
\int_{U\cap\C_j}G_L^{*s}(\varphi)(u+jv)\psi(u+jv)dudv=\int_{U\cap \C_j}\varphi(u+jv)G_L(\psi)(u+jv)dudv
$$
for every $\varphi\in\mathcal{SD}_R(U)$ and $\psi\in\mathcal{SD}_L(U)$. By the previous computations it is clear that the integrals on both sides do not dipend on $j$.

\begin{Dn}
The global (distributional) adjoints are given by
$$
G_L^*(\varphi)(q)=-G_R(\varphi)(q)-4\varphi(q)\underline{q}
$$
$$
G_R^*(\varphi)(q)=-G_L(\varphi)(q)-4\underline{q}\varphi(q)
$$
whereas the slice-adjoints are given by
$$
G_L^{*s}(\varphi)(q)=-G_R(\varphi)(q)-2\varphi(q)\underline{q}
$$
$$
G_R^{*s}(\varphi)(q)=-G_L(\varphi)(q)-2\underline{q}\varphi(q)\;.
$$
\end{Dn}

\begin{La} We have that
$$
G_L(S^{-1}_L(s,\cdot))=2\pi j|\underline{s}|^2\delta_s
$$
in the sense of distribution.
\end{La}
\begin{proof}Fix $s\in\hh$ and let $z_j=s_0+|\underline{s}|j$ (and $\overline{z}_j=s_0-|\underline{s}|j$); for $\varepsilon>0$, set
$$
W_{j,\varepsilon}=\{w\in \C_j\ :\ \min\{|w-z_j|, |w-\overline{z}_j|\}\geq\varepsilon\}
$$
and
$$
\Gamma_j=\{w\in \C_j\ :\ |w-z_j|=\varepsilon\},\qquad \overline{\Gamma}_j=\{w\in \C_j\ :\ |w-\overline{z}_j|=\varepsilon\}\;.
$$
For every $\varphi\in\mathcal{SD}_R(\hh)$, we have
\begin{eqnarray*}\int_{W_{j,\varepsilon}}G_L^{*s}(\varphi)(u+jv)S^{-1}_L(s,u+jv)dudv&=
&\int_{W_{j,\varepsilon}}\varphi(u+jv)G_L(S^{-1}_L)(s,u+jv)dw\wedge d\bar{w}\\
&&+\int_{\Gamma_j\cup\overline{\Gamma_j}}y^2\varphi(u+jv)S^{-1}_L(s,u+jv)(du+jdv)
\end{eqnarray*}
 As $S^{-1}_L(s,p)$ is a left slice regular function for $p\in W_{j,\varepsilon}$,
  the first integral vanishes; in order to compute the remaining one as $\varepsilon\to 0$, we note that, setting $j_s=\underline{s}|\underline{s}|^{-1}$,
$$S^{-1}_L(s,u+jv)=\dfrac{1}{2}\left(\dfrac{1}{(u+jv)-z_j}(1-jj_s)+\dfrac{1}{(u+jv)-\overline{z}_j}(1+jj_s)\right)$$
and
$$\lim_{\varepsilon\to 0}\int_{\Gamma_j}v^2\varphi(u+jv)\dfrac{1}{(u+jv)-z_j}(du+jdv)=2\pi j |\underline{s}|^2\phi(z_j)$$
$$\lim_{\varepsilon\to 0}\int_{\overline{\Gamma}_j}y^2\varphi(u+jv)\dfrac{1}{(u+jv)-\overline{z}_j}(du+jdv)=2\pi j |\underline{s}|^2\phi(\overline{z}_j)\;.$$
Hence
$$\int_{\Gamma_j\cup\overline{\Gamma_j}}y^2\varphi(x+jy)S^{-1}_L(s,u+jv)(du+jdv)=2\pi j |\underline{s}|^2\phi(s)\;.$$

This computation then implies that
$$G_L(S^{-1}_L(s,\cdot))=2\pi j|\underline{s}|^2\delta_s$$
in the sense of distributions.\end{proof}

We performed the computations on a slice $U\cap\C_j$; however, using the global adjoint, we could have carried on the same steps with respect to the Lebesgue measure on $U$, getting to the same result. We are now ready to solve the quaternionic global equation on a bounded domain.

\begin{Tm}\label{teo_debar_slice}Let $U$ be a bounded axially symmetric set, $\mathcal{X}$ a two sided quaternionic Banach space and  $V:U\to \mathcal{X}$ a left slice-$L^\infty$ function; define
$$f(p)=\dfrac{1}{2\pi j}\int_{U\cap\C_j}S_L^{-1}(s,p)V(s)ds\wedge d\bar{s}$$
for $p\in\H$.
Then $f(q)$ is a slice-$L^\infty$ function which solves $G_L(f)(q)=|\underline{q}|^2V(q)$ for almost every $q\in U$.\end{Tm}
\begin{proof} The $L^\infty$-norm of $f$ is bounded by
$$\|V\|_{\infty}\sup_{p\in U}\|S^{-1}_L(\cdot, p)\|_{L^1}$$
which can be easily computed similarly to the $L^1$ norm of the Cauchy kernel in one complex variable.

Moreover, the function $p\mapsto S^{-1}_L(s,p)V(s)$ is a left-slice function, for every $s\in U$, $[s]\neq[p]$; therefore, by the linearity of the integral, also the resulting function $f(p)$ is a left-slice function.

For every $\varphi\in\mathcal{SD}_R(U)$, we compute
\[
\begin{split}
\int_{U\cap\C_j}\varphi(p)G_Lf(p)d\lambda(p)&=\int_{U\cap\C_j}G_L^*\varphi(p)f(p)d\lambda(p)
\\
&
=\dfrac{1}{2\pi}\int_{U\cap\C_j}G^*_L\varphi(p)\int_{U\cap\C_j}S_L^{-1}(s,p)V(s)d\lambda(s)d\lambda(p)
\\
&
=\dfrac{1}{2\pi}\int_{U\cap\C_j}\left[\int_{U\cap\C_j}G^*L\varphi(p)S_L^{-1}(s,p)d\lambda(p)\right]V(s)d\lambda(s)
\\
&
=\int_{U\cap\C_j}\varphi(s)|\underline{s}|^2V(s)d\lambda(s)
\end{split}
\]
hence $G_Lf(p)\equiv |\underline{p}|^2V(p)$ almost everywhere on $U\cap\C_j$, hence on $U$.\end{proof}

As we noted in the beginning of the section, if $G_L(f)(p)=|\underline{p}|^2V(p)$, then
$$\overline{\partial}_jf(u+jv)=V(u+jv)$$ for almost every $(u+jv)\in U\cap \C_j$. Combining this result with Runge theorem, we obtain the solvability of the quaternionic global Cauchy-Riemann equation on any axially symmetric domain.

\begin{Tm}\label{teo_debar_open}
Given an axially symmetric slice-domain $U$ and a left slice-continuous function $V:U\to \mathcal{X}$, there exists a left slice-continuous function $f:U\to\mathcal{X}$ such that $G_L(f)(q)=|\underline{q}|^2V(q)$ for all  $q\in U$.
\end{Tm}
\begin{proof} Let $\{U_n\}$ an increasing sequence of axially symmetric domains such that $\overline{U}_n\subset\stackrel{\circ}{U}_{n+1}\Subset U$, $U_n$ has $\mathcal{C}^1$ boundary and each bounded connected component of $\C_j\setminus U_n$ (if any) contains at least one point of $\C_j\setminus U$, for each $n\in\N$ and each $j\in\S$. Moreover, let $\bigcup_n U_n=U$.

By Theorem \ref{teo_debar_slice}, we have a collection of functions $g_n:U_n\to\H$ such that
$G_L(g_n)(q)=|\underline{q}|^2V(q)$ for all $q\in U_n$.

We construct inductively another family such that $G_L(f_n)(q)=|\underline{q}|^2V(q)$
for all $q\in U_n$. Set $f_0=g_0$, $f_1=g_1$ and, if we have chosen $f_0,\ldots, f_k$,
we define $f_{k+1}$ as follows: the difference $f_k-g_{k+1}$ is left-slice regular on $U_k\supset\overline{U}_{k-1}$; hence, by
Runge approximation theorem (which is proved in \cite{RUNGE} for scalar functions, but can be easily generalized to functions with values in a Banach space) we have a left-slice regular function $h:U\to\mathcal{X}$ such that
$$\max_{q\in\overline{U}_{k-1}}|f_k(q)-g_{k+1}(q)-h(q)|\leq \dfrac{1}{2^{k+1}}\;.$$
We set $f_{k+1}=g_{k+1}+h\vert_{U_{n+1}}$.

Now, the sequence $f_k$ converges uniformly on compact sets on $U$ to a left slice-regular function $f:U\to\mathcal{X}$; moreover, also the sequence $\{f_{n}-f_k\}_{n\geq k}$ converges uniformly on $U_k$, to $f-f_k$ which then is slice-regular as well. So $G_L(f-f_k)=0$ on $U_k$, implying that
$$G_L(f)=G_L(f_k)=|\underline{q}|^2V$$
on $U_k$. \end{proof}

\subsection{Cocycles and Mittag-Leffler theorem}

Compare
Mittag-Leffler theorem in \cite{Sheaves}.

Slice hyperholomorphic functions admit a power series expansion in terms of suitable polynomials, see \cite{gps}:
\begin{Dn}\label{sphseries}
Let $q_0\in \mathbb H$. For any sequence $c_n\in\mathbb H$, $n\in\mathbb{Z}$, the series
\[
\sum_{n\in\mathbb Z} (q^2-2{\rm Re}(q_0) q+|q_0|^2)^n(c_{2n}+(q-q_0)c_{2n+1})
\]
is called the spherical Laurent series centered at $q_0$ associated with $\{c_n\}$, $n\in\mathbb{Z}$. If $c_n = 0$ for $n < 0$, then it is called the spherical series centered at $q_0$ associated with  $\{c_n\}$, $n\in\mathbb{Z}$.
\end{Dn}
To study the convergence of this series, one need the pseudo-distance
$$
d(q,p)=\sqrt{q^2-2{\rm Re}(q_0) q+|q_0|^2}.
$$
We then define the so-called Cassini ball with center at $q_0$ as
$$
U(q_0, r)=\{q\in\mathbb H\ :\ d(q,q_0)<r\},
$$
and the so-called Cassini shell
$$
U(q_0, r_1,r_2)=\{q\in\mathbb H\ :\ r_1<d(q,q_0)<r_2\}.
$$
\begin{Tm}\label{spherTm}
Let $\{c_n\}$, $n\in\mathbb{Z}$ and let
$$
r_1=\overline{\lim}_{n\to +\infty} \| c_{-n}\|^{1/n},\qquad
\frac{1}{r_2}=\overline{\lim}_{n\to +\infty} \| c_{n}\|^{1/n}.
$$
Consider the spherical Laurent series
\[
f(q)=\sum_{n\in\mathbb Z} (q^2-2{\rm Re}(q_0) q+|q_0|^2)^n(c_{2n}+(q-q_0)c_{2n+1}).
\]
If there exists $n<0$ such that $c_n\not=0$ then the domain of convergence is the Cassini shell
$U(q_0,r_1,r_2)$. If all $c_n=0$ for all $n<0$ then $f(q)$ is a spherical series and its domain of convergence is the Cassini ball $U(q_0,r_2)$. When the domain $U$ of convergence of $f$ is nonempty, then $f(q)$ is slice hyperholomorphic in $U$.
\end{Tm}

We also recall the definition of order of a function at a singular point or sphere, see \cite{gps}.
\begin{Dn}\label{SVILSFER}
Let $f:\, U\to\mathbb H$ be a slice hyperholomorphic function, let $q_0$ be a singularity for $f$ and let
\[
f(q)=\sum_{n\in\mathbb Z} (q^2-2{\rm Re}(q_0) q+|q_0|^2)^n(c_{2n}+(q-q_0)c_{2n+1}).
\]
be the spherical Laurent expansion of $f$ at $q_0$. We define the spherical order of $f$ at $[q_0]$, and we denote it by ${\rm ord}_f([q_0])$, as the smallest even number $n_0\in\mathbb{N}$ such that $c_n=0$ for all $n<-n_0$. If no such $n_0$ exists then we set ${\rm ord}_f([q_0])=+\infty$.
\end{Dn}
The following result is also taken from \cite{gps}, section 9.
\begin{Tm}\label{1.10.3}
Let $\tilde U$ be an axially symmetric open set in $\mathbb H$, $q_0\in\tilde U$ and $U=\tilde U\setminus \{[q_0]\}$. Let $f:\, U\to\mathbb H$ be a slice hyperholomorphic function. Then to say that every point of $[q_0]$ is a removable singularity for $f$ is equivalent to
\begin{enumerate}
\item ${\rm ord}_f([q_0])=0$;
\item ${\rm ord}_f(q_0)={\rm ord}_f(\bar q_0)=0$;
\item There exists a neighborhood $U_{q_0}$ of $q_0$ in $U$ such that $\|f\|$ is bounded in $U_{q_0}\setminus \{[q_0]\}$.
\end{enumerate}
\end{Tm}

As an application of the solvability results for the Cauchy-Riemann equation, we translate to the quaternionic setting the traditional proof of Mittag-Leffler theorem, taking also the occasion to state it for functions with values in a Banach spaces.

\begin{Dn}Let $U\subseteq H$ be an axially symmetric domain and $\{U_k\}_{k\in K}$ an open covering of $U$ made of axially symmetric sets. A family $\{f_{hk}\}_{h,k\in K}$ of slice regular functions $f_{hk}:U_h\cap U_k\to\mathcal{X}$ is called a \emph{cocycle} if
$$f_{hk}+f_{kl}=f_{hl}$$
on $U_h\cap U_k\cap U_l$ whenever this set is non-empty.\end{Dn}

\begin{La}\label{Cousin}Let $U$ and $\{U_k\}_{k\in K}$ be as above; given a cocycle $\{f_{hk}\}_{h,k\in K}$, there exists a family $\{g_h\}_{h\in K}$ of slice-regular functions $g_h:U_h\to\mathcal{X}$ such that
$$f_{hk}=g_h-g_k$$
on $U_h\cap U_k$ whenever this set is non-empty.\end{La}
\begin{proof}Following the classical proof, we consider a (real-valued) partition of unity $\{\chi_h\}_{h\in K}$, where each $\chi_h$ is a slice function, subordinated to the open covering $\{U_h\}_{h\in K}$ and define slice-smooth maps $\phi_h:U_h\to \mathcal{X}$ by setting
$$\phi_h=-\sum_{k\in K}\chi_kf_{kh}\;.$$
Then $\phi_h-\phi_k=f_{hk}$ on $U_h\cap U_k$, which implies that $G_L(\phi_h-\phi_k)=0$ on $U_h\cap U_k$; hence, the functions $G_L(\phi_h)$ glue into a global function $|\underline{q}|^2\psi(q)$. The function $\psi$ is actually slice-smooth: the $\phi_h$'s are slice-functions, so the functions $\overline{\partial}_j\phi_h$ (in principle defined only on $U_h\cap\C_j$ glue into a slice-smooth function on $U_h$, which is $|\underline{q}|^{-2}G_L(\phi_h)(q)$.

By Theorem \ref{teo_debar_open}, there exists a function $u:U\to\mathcal{X}$ such that $G_L(u)(q)=|\underline{q}|^2\psi(q)$. Setting $g_h=\phi_h-u$ and recalling Remark \ref{rmk1}, we complete the proof.\end{proof}

\begin{Tm}\label{MittagLeffler}
Let $U$ be an axially symmetric domain, $Z$ a discrete and closed subset of $U$, and assume that, for each $p\in Z$, a slice-regular function $f_w:\hh\setminus [w]\to\mathcal{X}$ of the form
as in Definition \ref{SVILSFER}
is given. Then there exists a slice-regular function $f:U\setminus[Z]\to\mathcal{X}$ such that, for each $w\in Z$, $f_w$ is the principal part of the Laurent expansion of $f$ at $w$.
\end{Tm}
\begin{proof} We can find a collection of axially symmetric open sets $\{U_w\}_{w\in Z}$ such that $U_w$ is a neighbourhood of $[w]$ and $U_w\cap U_{w'}=\emptyset$ if $w\neq w'$. Set $U_1=\bigcup_{w\in Z} U_w$ and $U_2=U\setminus[Z]$.
Setting $g(q)=f_w(q)$ for $q\in U_w\setminus [w]$, we define a slice-regular function on $U_1\cap U_2$. By Lemma \ref{Cousin}, there exist $h_i:U_i\to\mathcal{X}$ such that $g=h_1-h_2$.

It is easy to verify that $f=-h_2$ has the required properties.\end{proof}

\end{document}